\documentclass[10pt]{article}
\usepackage{graphicx}
\newtheorem{theorem}{Theorem}
\newtheorem{proposition}{Proposition}
\newtheorem{lemma}{Lemma}
\newtheorem{corollary}{Corollary}

\input tcilatex

\begin{document}

\begin{titlepage}

\vskip0.5truecm

\begin{center}

{\LARGE \bf On $C^r$-generic twist maps of ${\rm T^2}$}

\end{center}

\vskip  0.4truecm

\centerline {{\large Salvador Addas-Zanata}}

\vskip 0.2truecm

\centerline { {\sl Instituto de Matem\'atica e Estat\'\i stica }}
\centerline {{\sl Universidade de S\~ao Paulo}}
\centerline {{\sl Rua do Mat\~ao 1010, Cidade Universit\'aria,}} 
\centerline {{\sl 05508-090 S\~ao Paulo, SP, Brazil}}
 
\vskip 0.7truecm

\begin{abstract}

We consider twist diffeomorphisms of the torus, 
$f:{\rm T^2\rightarrow T^2,}$ and their vertical rotation intervals 
$\rho _V(\widehat{f})=[\rho _V^{-},\rho _V^{+}],$ where $\widehat{f}$ is a
lift of $f$ to the vertical annulus or cylinder. We show that $C^r$-generically 
for any $r\geq 1$, both extremes of the rotation interval are rational and locally constant 
under $C^0$-perturbations of the map. Moreover, when $f$ is area-preserving, $C^r$-generically 
$\rho _V^{-}<\rho _V^{+}.$ 

Also, for any twist map $f$, $\widehat{f}$ a lift of $f$ to the cylinder,  
if $\rho _V^{-}<\rho _V^{+}=p/q$, then there are 
two possibilities: either $\widehat{f}^q(\bullet)-(0,p)$ maps a simple essential loop into the connected component of its 
complement which is below the loop, or it satisfies the Curve Intersection Property. 
In the first case, $\rho _V^{+} \leq p/q$ in a $C^0$-neighborhood of $f,$ and in the second case,
we show that $\rho _V^{+}(\widehat{f}+(0,t))>p/q$ for all $t>0$ 
(that is, the rotation interval is ready to grow).  
Finally, in the $C^r$-generic case, assuming that $\rho _V^{-}<\rho _V^{+}=p/q,$ 
we present some consequences of the existence of the free loop for $\widehat{f}^q(\bullet)-(0,p)$,
related to the description and shape of the attractor-reppeler pair that exists in the annulus. 
The case of a $C^r$-generic transitive twist diffeomorphism (if such a thing exists) 
is also investigated.

\end{abstract} 

\vskip 0.3truecm

\vskip 2.0truecm


\hskip1.8truecm full and partial meshes, Pixton's Theorem

\vskip 0.8truecm

\noindent{\bf e-mail:} sazanata@ime.usp.br 

\vfill
\hrule
\noindent{\footnotesize{The author would like to thank FAPESP (project 
2023/07076-4) for partial funding.}}

\end{titlepage}

\baselineskip=6.2mm

\section{Introduction and main results}

Understanding the relationship between the dynamics of homeomorphisms of
closed oriented surfaces in the identity isotopy class, and their
(homological) rotation sets has been a very active field of research.
Several results concerning these relations have been proved, probably
starting with \cite{mz}, \cite{llibremack} and \cite{franks1989}.

The aim of the present paper is to study this problem in the torus, mostly
for another homotopy class, the so called Dehn twists. In this case, the
rotation set is only one dimensional (see \cite{alunofranks} and \cite{non04}%
), and similarly to the isotopic to the identity class, it is convex. So,
the only possibilities are: a point or a non-degenerate closed interval. We
will deal mostly with the second possibility. In the majority of our
results, we assume an additional dynamical property, the twist condition.
This means that, if we consider a lift of the torus map to the plane, 
then the image of any vertical line projects injectively onto the horizontal
coordinate.

This class of mappings has very rich dynamics, as it contains the well-know
standard mapping, $S_M:{\rm T^2\rightarrow T^2,}$ which writes in flat
coordinates, as${\rm :}$ 
\begin{equation}
\label{stmap}S_M:\left\{ 
\begin{array}{l}
x^{\prime }=x+y+\frac k{2\pi }\sin (2\pi x) 
\text{ mod 1} \\ y^{\prime }=y+\frac k{2\pi }\sin (2\pi x )\text{ mod 1}%
 
\end{array}
\right. 
\end{equation}

More precisely, our aim is to study the following problems:

\begin{itemize}
\item  for a $C^r$-generic (for any $r\geq 1$) twist diffeomorphism of ${\rm %
T^2,}$ what can we say about its rotation set? Is it locally constant (does
mode-locking happen \cite{calvezeu})? We also consider the same questions in
the area-preserving world;

\item  find some property, that, when satisfied, implies that a general
twist diffeomorphism of {\rm T$^2$} has a locally constant rotation set, and
when it is not satisfied, show that the rotation set can be changed by
arbitrarily small perturbations, in any differentiability;

Using the results obtained related to the above questions, we further study
some dynamical consequences of locally constant extremes for the rotation
set, in $C^r$-genericity. We also present some consequences of non-existence
of periodic open disks, that is of transitivity.
\end{itemize}

\begin{description}
\item[Notations and definitions]{\bf:}

\vskip0.2truecm

1) Let $(x,y)$ denote coordinates in the flat torus ${\rm T^2=I%
\negthinspace  R^2/Z\negthinspace  \negthinspace  Z^2,}$ $(\widehat{x},
\widehat{y})$ in the annulus ${\rm T^1\times I\negthinspace  R,}$ and let $(
\widetilde{x},\widetilde{y})$ denote coordinates in ${\rm I\negthinspace  %
R^2.}$ Let $p_1,p_2$ from ${\rm T^2}$ or ${\rm T^1\times I\negthinspace  R}$
or ${\rm I\negthinspace  R^2},$ to ${\rm I\negthinspace  R}$ be the standard
projections, respectively in the horizontal and vertical coordinates, and
let $\pi :{\rm I\negthinspace  R^2\rightarrow }$ ${\rm T^1\times 
I\negthinspace  R},$ $\tau :{\rm T^1\times \rm I\negthinspace  R\rightarrow }$ ${\rm T^2}$ 
and $p:{\rm I\negthinspace  R^2\rightarrow T^2}$ be the
covering mappings.

2) Define ${\rm Diff_k^r(I\negthinspace  R^2)=\{{\it \widetilde{f}}:I\negthinspace  %
R^2\rightarrow I\negthinspace  R^2:}$ $\widetilde{f}$ is a $C^r$%
-diffeomorphism of the plane ($r\geq 0$) such that for any pair of integers $%
n,m,$ $\widetilde{f}(\widetilde{x}+n,\widetilde{y}+m)=\widetilde{f}(
\widetilde{x},\widetilde{y})+(n+km,m),$ for all $(\widetilde{x},\widetilde{y}%
)\in {\rm I\negthinspace  R^2,}$ where $k$ is an integer$\}$.

3) Similarly, let ${\rm Diff_k^r(T^1\times I\negthinspace  R)=\{{\it \widehat{f}}:T^1\times
I\negthinspace  R\rightarrow T^1\times I\negthinspace  R:}$ $\widehat{f}$ is a 

\noindent $C^r$-diffeomorphism ($r\geq 0$) lifted by some element $\widetilde{f}\in 
{\rm Diff_k^r(I\negthinspace  R^2)\}}$ and, ${\rm Diff_k^r(T^2)=}$
${\rm \{{\it f}:T^2\rightarrow T^2:}$ $f$ is a $C^r$-diffeomorphism ($r\geq
0 $) lifted by some element $\widetilde{f}\in {\rm Diff_k^r(I\negthinspace  %
R^2)\}}.$ The union of ${\rm Diff_k^r(T^2)}$ for all $k\neq 0$ consists of the subset 
of torus $C^r$-diffeomorphisms homotopic to Dehn twists, and 
${\rm Diff_0^r(T^2)}$ is the set of torus $C^r$-diffeomorphisms homotopic to the identity.

4) For any $\widetilde{f}\in {\rm Diff_k^r(I\negthinspace  R^2),}$ $r\geq 1$
and $k\neq 0,$ if $k\times \partial _{\widetilde{y}}$ $(p_1\circ \widetilde{f}%
(\widetilde{x},\widetilde{y}))>0$ everywhere, then we say that $\widetilde{f}
$ is the lift of a torus (or annulus) twist map. Twist is to the right if $%
k>0$ and to the left if $k<0.$ When $k\neq 0$ and only twist maps are
considered in the previous subsets of diffeomorphisms, we add the letter $t:%
{\rm Diff_{t,k}^r(I\negthinspace  R^2),}$ ${\rm Diff_{t,k}^r(T^1\times I%
\negthinspace  R)}$ and ${\rm Diff_{t,k}^r(T^2).}$

5) For any $\widehat{f}\in {\rm Diff_k^0(T^1\times I\negthinspace  R),}$ $%
k\neq 0,$ and $z\in {\rm T^2,}$ we define the vertical rotation number of $%
z$ as (when the limit exists):

\item  
\begin{equation}
\label{vertrot}\rho _V(z)=\stackunder{n\rightarrow \infty }{\lim }\frac{%
p_2\circ \widehat{f}^n(\widehat{z})-p_2(\widehat{z})}n,\text{ for any } 
\widehat{z}\in \tau ^{-1}(z) 
\end{equation}

and the vertical rotation interval of $\widehat{f}$ is given by: 
\begin{equation}
\label{mz}\rho _V(\widehat{f})=\stackrel{\infty }{\stackunder{i=1}{\cap }} 
\overline{\stackunder{n\geq i}{\cup }\left\{ \frac{p_2\circ \widehat{f}^n( 
\widehat{z})-p_2(\widehat{z})}n:\widehat{z}\in {\rm T^1\times I\negthinspace  %
R}\right\} }, 
\end{equation}
which is a closed interval, possibly a point. Initially, in \cite{alunofranks}
and \cite{non04}, it was proved that all rational points $p/q$ in the
interior of $\rho _V(\widehat{f})$ are realized by $q$-periodic orbits in
the torus, similarly to the identity isotopy class. Moreover, in \cite
{alunofranksmisiu}, the function ${\rm Diff_k^0(T^1\times I\negthinspace  R)}%
\ni \widehat{f}\mapsto \rho _V(\widehat{f})=[\rho _V^{-},\rho _V^{+}]$
was shown to vary continuously (see also \cite{qualita}).

6) For any integer $k$ and $\widehat{f}\in {\rm Diff_k^0(T^1\times I%
\negthinspace  R),}$ we say that $\widehat{f}$ satisfies the Curve Intersection
Property (C.I.P.), if $\widehat{f}(\gamma )\cap \gamma \neq \emptyset $ for
all homotopically non-trivial simple closed curves $\gamma \subset {\rm T^1\times I\negthinspace  R.}$

7) Let $f:X\rightarrow X$ be a map from a metric space $(X,d),$ $n>0$ be an
integer and $\varepsilon >0$ a fixed real number. A sequence 
$\{x_0,x_1,...,x_n\}\subset X$ is an $\varepsilon $-pseudo orbit if 
$d(f(x_i),x_{i+1})<\varepsilon $ for all $i=0,1,...,n-1.$

8) We say that $D$, an open subset of ${\rm T^2}$ is essential if $D$ contains 
an homotopically non-trivial simple closed curve in the torus. We say that 
$D$ is fully essential 
if it contains two homotopically non-trivial simple closed curves which are not in
the same homotopy class. In this case, $D$ contains closed curves in all homotopy 
classes and $D^c$ is contained in a disjoint union of open disks. Moreover, we say that 
$D$ is inessential if it is not essential.

9) Given a closed subset $K$ of ${\rm T^2}$, ${\rm Filled}(K)$ is given by the union 
of $K$ with all the connected components of its complement which are inessential.

\end{description}

Now we are ready to state our main results. The first one explains when the 
vertical rotation interval can, or cannot change under perturbations.

\begin{theorem}
\label{pic} Let $k\neq 0$ and $r\geq 1$ be integers and $\widehat{f}\in {\rm %
Diff_{t,k}^r(T^1\times I\negthinspace  R)}.$ The following assertions are 
equivalent: 

\begin{itemize}

\item $\rho _V^{+}:{\rm Diff_{t,k}^r(T^1\times I\negthinspace  R)\rightarrow I%
\negthinspace  R}$ has a local maximum at $\widehat{f};$

\item $\rho _V^{+}(\widehat{f})$ is equal to some rational number $p/q$ and $\widehat{f}%
^q(\bullet )-(0,p)$ is an annulus diffeomorphism which maps a homotopically
non-trivial simple closed curve $\gamma $ into $\gamma $$^{-},$ the
connected component of $\gamma ^c$ which is below $\gamma ;$ 
\end{itemize}

So, if $\rho _V^{+}:{\rm Diff_{t,k}^r(T^1\times I\negthinspace  %
R)\rightarrow I\negthinspace  R}$ has a local maximum at $\widehat{f},$ it
also has a local maximum at $\widehat{f}$ in ${\rm Diff_k^0(T^1\times I\negthinspace  R).%
}$
\end{theorem}

\begin{corollary}
\label{remark} Let $k\neq 0$ and $r\geq 1$ be integers and $\widehat{f}\in {\rm %
Diff_{t,k}^r(T^1\times I\negthinspace  R)}.$ If $\rho
_V^{-}(\widehat{f})<\rho _V^{+}(\widehat{f})=p/q,$ then the following equivalences hold:

\begin{itemize}

\item $\rho _V^{+}$ does not have a local maximum at $\widehat{f};$ 

\item $\widehat{f}^q(\bullet )-(0,p)$ satisfies the Curve
Intersection Property (C.I.P.); 

\item $\rho _V^{+}(\widehat{f}+(0,t))>p/q$ for all $t>0;$

\end{itemize}

Moreover, when C.I.P. holds, then $\widehat{f}^q(\bullet )-(0,p)$ has
periodic orbits of all rational rotation numbers in the annulus.
\end{corollary}

Of course, analogous versions related to $\rho _V^{-}(\widehat{f})$ also
hold.

The next lemma is an important tool in the proof of the above results. It is 
contained in Corollary 3 of \cite{euerg1} and Lemma 2 and Theorem 10 of \cite{qualita}.
The ideas behind its proof rely on the topological theory for twist maps developed by Le
Calvez in \cite{lecalvezcompte}.

\begin{lemma}
\label{euerg1lemma} Let $k\neq 0$ and $r\geq 1$ be integers and $\widehat{f}\in {\rm %
Diff_{t,k}^r(T^1\times I\negthinspace  R)}.$ The
function $s \mapsto \rho _V^{+}(\widehat{f}+(0,s ))$ is continuous
and non-decreasing. Moreover, if there exist maps $\widehat{h}\in {\rm %
Diff_{t,k}^1(T^1\times I\negthinspace  R)}$ arbitrarily $C^0$-close to $
\widehat{f}$ such that $\rho _V^{+}(\widehat{h})>\rho _V^{+}(\widehat{f}),$
then for all $s>0,$ $\rho _V^{+}(\widehat{f}+(0,s))>\rho _V^{+}( 
\widehat{f}).$ Analogously, if there exist maps $\widehat{h}\in {\rm %
Diff_{t,k}^1(T^1\times I\negthinspace  R)}$ arbitrarily $C^0$-close to $
\widehat{f}$ such that $\rho _V^{+}(\widehat{h})<\rho _V^{+}(\widehat{f}),$
then for all $s<0,$ $\rho _V^{+}(\widehat{f}+(0,s))<\rho _V^{+}( 
\widehat{f}).$
\end{lemma}

{\bf Remark: }As above, an analogous statement holds for $\rho _V^{-}.$

In other words, given a twist diffeomorphism of the torus, for which one
extreme of the vertical rotation interval is not locally constant, then it
varies under arbitrarily small vertical translations of the diffeomorphism.
Or equivalently, if $s\mapsto \rho _V^{+(-)}(\widehat{f}+(0,s))$
is locally constant in a neighborhood of zero, then $\rho _V^{+(-)}$ is
locally constant in a $C^1$-neighborhood of $\widehat{f},$ which implies by
Theorem \ref{pic} that $\rho _V^{+(-)}(\widehat{f})$ is some rational
number $p/q$ and the annulus diffeomorphism $\widehat{f}^q(\bullet )-(0,p)$
has a homotopically non trivial free simple curve $\gamma,$ that is, 
the image of $\gamma$ under $\widehat{f}^q(\bullet )-(0,p)$ is disjoint from $\gamma.$

Corollary \ref{remark} also says that, whenever for an one-parameter family of
twist maps of the torus, $(f_t)_{t\in I},$ with $\rho _V^{-}( 
\widehat{f}_t)<\rho _V^{+}(\widehat{f}_t)$ for all $t\in I,$ if 
$\rho _V^{+}(\widehat{f}_{t^{*}})$ is a rational number $p/q$ and for all 
$t>t^{*},$ sufficiently close to $t^{*},$ $\rho _V^{+}(\widehat{f}_t)>\rho
_V^{+}(\widehat{f}_{t^{*}}),$ then $(\widehat{f}_{t^{*}})^q(\bullet
)-(0,p) $ satisfies the C.I.P.

Then next two theorems describe vertical rotation intervals for $C^r$-generic 
twist maps of the torus.

\begin{theorem}
\label{conjrot} Let $k\neq 0$ and $r\geq 1$ be integers. The set $O{\rm _{t,k}^r}({\rm %
T^1\times I\negthinspace  R})$ of maps $\widehat{f}\in {\rm %
Diff_{t,k}^r(T^1\times I\negthinspace  R)}$ such that $\rho _V$ is constant
in a neighborhood of $\widehat{f},$ is open and dense. Furthermore, the
endpoints of $\rho _V(\widehat{f})$ are rational and $\rho _V$ is constant
in a neighborhood of $\widehat{f}\in {\rm Diff_k^0(T^1\times I\negthinspace  %
R).}$
\end{theorem}

In other words, for twist maps, $C^r$-generically for any $r\geq 1,$ the vertical rotation
interval is locally constant, and its extremes are rational numbers. There
is also an area-preserving version of the above theorem. Note that 
${\rm Diff_{t,k,Leb}^r(T^1\times I\negthinspace  R)\stackrel{def.}{=}\{}$%
area-preserving elements of ${\rm Diff_{t,k}^r(T^1\times I\negthinspace  %
R)\}}$:

\begin{theorem}
\label{conjrotarea} Let $k\neq 0$ and $r\geq 1$ be integers. The set $O{\rm %
_{t,k,Leb}^r}({\rm T^1\times I\negthinspace  R})$ of maps $\widehat{f}\in 
{\rm Diff_{t,k,Leb}^r(T^1\times I\negthinspace  R)}$ 
such that $\rho _V$ is constant in a neighborhood of $\widehat{f},$
is open and dense. Furthermore, the endpoints of $\rho _V(\widehat{f})$ are
different rational numbers and $\rho _V$ is constant in a neighborhood of $
\widehat{f}\in {\rm Diff_k^0(T^1\times I\negthinspace  R).}$
\end{theorem}

The only difference between the general setting and the area-preserving, is
that in the last one, $C^r$-generically, vertical rotation intervals have
interior.

Theorems 1,2 and 3 are the foundation on which our next results rely: a
description of the attractor-reppeller pair in the annulus, that always exists in
the $C^r$-generic case. 


Let $k\neq 0$ and $r\geq 1$ be integers and $\widehat{f}\in O{\rm _{t,k}^r}({\rm %
T^1\times I\negthinspace  R})\cap \chi ^r({\rm T^2})$ or in the area-preserving case,
$\widehat{f}\in O{\rm _{t,k,Leb}^r}({\rm T^1\times I%
\negthinspace  R}) \cap \chi _{{\rm Leb}}^r({\rm T^2}),$ 
the open and dense sets from Theorems \ref{conjrot} and 
\ref{conjrotarea} intersected with the residual sets 
$\chi ^r({\rm T^2})$ and $\chi _{{\rm Leb}}^r({\rm T^2})$ (from Theorem \ref{comeconovo}), 
which are contained in general Kupka-Smale, or in
area-preserving Kupka-Smale $C^r$-diffeomorphisms respectively, for which
closures of stable and unstable branches of hyperbolic periodic saddles and
homoclinic and heteroclinic intersections vary continuously with
perturbations of the diffeomorphisms. See \cite{pixton}, pages 370-372 for
more information.
  
So $\rho _V(\widehat{f})=[r/s,p/q]$ for rational numbers $%
r/s\leq $ $p/q$ and $\rho _V$ is locally constant in a neighborhood of $
\widehat{f}.$ From now on, we assume that $r/s<p/q$ (this is always the case
for area-preserving generic twist diffeomorphisms). 

Theorem \ref{pic} says that  
$\widehat{f}^q(\bullet )-(0,p)$ has a free homotopically non-trivial simple
closed curve $\gamma _{p/q}\subset {\rm T^1\times I\negthinspace  R}$, such
that: 
$$
\widehat{f}^q(\gamma _{p/q})-(0,p)\subset \gamma _{p/q}^{-} 
$$

Something which implies the existence of an attractor-repeller pair for $
\widehat{f}^q(\bullet )-(0,p).$ The attractor $A_{p/q}$ is contained in $%
\gamma _{p/q}^{-}$ and the repeller $R_{p/q}$ is contained in $\gamma
_{p/q}^{+}.$ In the next result, we describe this pair ($f\in {\rm %
Diff_{t,k}^r(T^2)}$ is the torus map lifted by $\widehat{f}$).

\begin{theorem}
\label{locconst} Under the previous hypotheses, there exists a hyperbolic
$f$-periodic saddle $z_{p/q}\in {\rm T^2}$ of vertical rotation number $p/q$
such that if $\widehat{z}_{p/q}\in {\rm T^1\times I\negthinspace  R}$ is any
lift of $z_{p/q}$ to the annulus, then $W^u(\widehat{z}_{p/q})$ is bounded
from above as a subset of the annulus and unbounded from below, $W^s( 
\widehat{z}_{p/q})$ is unbounded from above and bounded from below.
Moreover, $W^u(\widehat{z}_{p/q})$ has a transversal intersection with $W^s( 
\widehat{z}_{p/q}-(0,1)),$ and if $\widehat{z}_{p/q}$ and $\widehat{z}%
_{p/q}-(0,1)$ are both above $\gamma _{p/q},$ and $\widehat{z}_{p/q}-(0,2)$
is below, then $A_{p/q}$ is contained in $\overline{W^u(\widehat{z}_{p/q})}$ 
$\cup (\overline{W^u(\widehat{z}_{p/q})}^{b.above}),$ where the last set is
the union of all (open) connected components of $(\overline{W^u(\widehat{z}%
_{p/q})})^c$ which are bounded from above. Moreover, $A_{p/q}\supseteq
\left[ \overline{W^u(\widehat{z}_{p/q})}\cup (\overline{W^u(\widehat{z}%
_{p/q})}^{b.above})\right] $$-(0,2).$ Similarly, if $\widehat{z}_{p/q}$ and $
\widehat{z}_{p/q}+(0,1)$ are both below $\gamma _{p/q},$ and $\widehat{z}%
_{p/q}+(0,2)$ is above, then $R_{p/q}$ is contained in $\overline{W^s( 
\widehat{z}_{p/q})}$ $\cup (\overline{W^s(\widehat{z}_{p/q})}^{b.below}),$
where the last set is the union of all connected components of $(\overline{%
W^s(\widehat{z}_{p/q})})^c$ which are bounded from below, and $%
R_{p/q}\supseteq \left[ \overline{W^s(\widehat{z}_{p/q})}\cup (\overline{%
W^s( \widehat{z}_{p/q})}^{b.below})\right] $$+(0,2)$.
\end{theorem}

Our main interest in the previous result is to apply it in case $f:{\rm %
T^2\rightarrow T^2}$ is transitive. Ideally, we would like to understand if
such a $C^r$-generic twist diffeomorphism (for any $r\geq 1$), mostly in the
area-preserving case, could be transitive. Our hope is that for
large $r,$ it cannot. And this final result can be seen as an attempt to
start a list of consequences of generic transitivity, that ultimately, would
lead to a contradiction.

\begin{corollary}
\label{locconsttran} Still under the previous theorem's hypotheses, if we
assume $f$ to be transitive, then the following improvement holds:

\begin{itemize}
\item  if $\widehat{z}_{p/q}$ and $\widehat{z}_{p/q}-(0,1)$ are both above $%
\gamma _{p/q},$ and $\widehat{z}_{p/q}-(0,2)$ is below, then $A_{p/q}$ is
contained in $\overline{W^u(\widehat{z}_{p/q})}$ and contains $\overline{%
W^u( \widehat{z}_{p/q})}-(0,2);$

\item  similarly, if $\widehat{z}_{p/q}$ and $\widehat{z}_{p/q}+(0,1)$ are
both below $\gamma _{p/q},$ and $\widehat{z}_{p/q}+(0,2)$ is above, then $%
R_{p/q}$ is contained in $\overline{W^s(\widehat{z}_{p/q})}$ and contains $
\overline{W^s(\widehat{z}_{p/q})}+(0,2);$
\end{itemize}

Moreover, both $\overline{W^u(\widehat{z}_{p/q})}$ and $\overline{W^s( 
\widehat{z}_{p/q})}$ have no interior points and their complements are
connected, although $\overline{W^u(z_{p/q})}=\overline{W^s(z_{p/q})}={\rm T^2}$.
\end{corollary}

\section{Basic tools}

\subsection{Some results for twist maps}

\subsubsection{Le Calvez's topological theory}

The results below 
can be found in \cite{lecalvez1} and \cite{lecalvez2}. Let $k\neq 0$ be an 
integer, $\widehat{f}\in {\rm Diff_{t,k}^1}(\rm T^1\times {\rm I}\negthinspace 
{\rm R})$ and $\widetilde{f}\in {\rm Diff_{t,k}^1}({\rm I}\negthinspace {\rm %
R^2)}$ be one of its lifts. For every pair of integers $(s,q),$ $q>0,$ we
define the following sets:

\begin{equation}
\label{Kpq} 
\begin{array}{c}
K 
{\rm _{lift}}(s,q)=\left\{ (\widetilde{x},\widetilde{y})\in {\rm I}%
\negthinspace {\rm R^2}\text{: }p_1\circ \widetilde{f}^q(\widetilde{x}, 
\widetilde{y})=\widetilde{x}+s\right\} \\ \text{ and } \\ K(s,q)=\pi \circ K%
{\rm _{lift}}(s,q) 
\end{array}
\end{equation}

Then we have the following:

\begin{lemma}
\label{Lcal1} The set $K(s,q)$ is compact and it has a unique connected component $%
C(s,q)$ that separates the ends of the annulus.
\end{lemma}

Next, we define the following functions on ${\rm T^1}$: 
$$
\begin{array}{c}
\mu ^{-}( 
\widehat{x})=\min \{p_2(\widehat{z})\text{: }\widehat{z}\in K(s,q)\text{ and 
}p_1(\widehat{z})=\widehat{x}\} \\ \mu ^{+}(\widehat{x})=\max \{p_2(\widehat{%
z})\text{: }\widehat{z}\in K(s,q)\text{ and }p_1(\widehat{z})=\widehat{x}\} 
\end{array}
$$
And similar functions for $\widehat{f}^q(K(s,q))$: 
$$
\begin{array}{c}
\nu ^{-}( 
\widehat{x})=\min \{p_2(\widehat{z})\text{: }\widehat{z}\in \widehat{f}%
^q\circ K(s,q)\text{ and }p_1(\widehat{z})=\widehat{x}\} \\ \nu ^{+}( 
\widehat{x})=\max \{p_2(\widehat{z})\text{: }\widehat{z}\in \widehat{f}%
^q\circ K(s,q)\text{ and }p_1(\widehat{z})=\widehat{x}\} 
\end{array}
$$

The following lemmas are very important in this theory:

\begin{lemma}
\label{graphu} Defining ${\rm Graph}$\{$\mu ^{\pm }$\}=\{$(\widehat{x},\mu
^{\pm }(\widehat{x})):\widehat{x}\in {\rm T^1}$\} we have: 
$$
{\rm Graph}\{\mu ^{-}\}\cup {\rm Graph}\{\mu ^{+}\}\subset C(s,q) 
$$
\end{lemma}

\begin{lemma}
\label{ofpre} The following equalities hold for all $\widehat{x}\in S^1:$ 


$\widehat{f}^q(\widehat{x},\mu ^{-}(\widehat{x}))=(\widehat{x}%
,\nu ^{+}(\widehat{x}))$ and $\widehat{f}^q(\widehat{x},\mu ^{+}(\widehat{x}%
))=(\widehat{x},\nu ^{-}(\widehat{x})).$
\end{lemma}


Now we recall ideas and results from \cite{lecalvezcompte}. Fix some $
\widetilde{f}\in {\rm Diff_{t,k}^1}({\rm I}\negthinspace {\rm R^2)}$ and let 
$\widehat{f}$ be the annulus map lifted by $\widetilde{f}$.

Given a triplet of integers $(s,p,q)$ with $q>0,$ if there is no point $( 
\widetilde{x},\widetilde{y})\in {\rm I}\negthinspace 
{\rm R^2}$ such that $\widetilde{f}^q(\widetilde{x},\widetilde{y})=( 
\widetilde{x}+s,\widetilde{y}+p),$ it can be proved that the sets $\widehat{f%
}^q\circ K(s,q)$ and $K(s,q)+(0,p)$ can be separated by the graph of a
continuous function $\sigma :{\rm T^1}\rightarrow {\rm I}\negthinspace {\rm R}$,
essentially because from all the previous results, either one of the
following inequalities must hold: 
\begin{equation}
\label{pos}\nu ^{-}(\widehat{x})-\mu ^{+}(\widehat{x})>p 
\end{equation}
\begin{equation}
\label{neg}\nu ^{+}(\widehat{x})-\mu ^{-}(\widehat{x})<p 
\end{equation}
for all $\widehat{x}\in {\rm T^1},$ where $\nu ^{+},\nu ^{-},\mu ^{+},\mu ^{-}$
were defined above.

Following Le\ Calvez \cite{lecalvezcompte}, we say that the triplet $(s,p,q)$
is positive (resp. negative) for $\widetilde{f}$ if $\widehat{f}^q\circ
K(s,q)$ is above (\ref{pos}) (resp. below (\ref{neg})) the graph of $\sigma
. $

Recall that: 
\begin{equation}
\label{Tgglin}\widetilde{f}(\widetilde{x},\widetilde{y})=(\widetilde{x}%
^{\prime },\widetilde{y}^{\prime })\Leftrightarrow \widetilde{y}=g( 
\widetilde{x},\widetilde{x}^{\prime })\text{ and }\widetilde{y}^{\prime
}=g^{\prime }(\widetilde{x},\widetilde{x}^{\prime }), 
\end{equation}
where $g,g^{\prime }$ are mappings from ${\rm I}%
\negthinspace 
{\rm R^2}$ to ${\rm I}\negthinspace {\rm R},$ 
which satisfy $g^{\prime }(\widetilde{x},\widetilde{x}^{\prime })=p_2\circ 
\widetilde{f}(\widetilde{x},g(\widetilde{x},\widetilde{x}^{\prime })).$

Before stating the next proposition, we need some extra definitions and 
simple facts.

If we define $\widetilde{h}_t(\widetilde{x},\widetilde{y})=(\widetilde{x}, 
\widetilde{y}+t),$ it is easy to see that for all $t\in {\rm I}%
\negthinspace 
{\rm R,}$ $\widetilde{h}_t(\widetilde{x},\widetilde{y})$ conjugates $( 
\widetilde{x},\widetilde{y})\mapsto \widetilde{f}(\widetilde{x}, 
\widetilde{y}+t)+(0,t)$ with $\widetilde{f}_t(\widetilde{x},\widetilde{y})%
\stackrel{def.}{=}\widetilde{f}(\widetilde{x},\widetilde{y})+(0,2t).$ And if
we denote as $g_t(\widetilde{x},\widetilde{x}^{\prime })$ and $g_t^{\prime
}( \widetilde{x},\widetilde{x}^{\prime })$ the mappings associated to $
\widetilde{f}(\widetilde{x},\widetilde{y}+t)+(0,t)$ in the way of (\ref
{Tgglin}), then $g_t(\widetilde{x},\widetilde{x}^{\prime })=g(\widetilde{x}, 
\widetilde{x}^{\prime })-t$ and $g_t^{\prime }(\widetilde{x},\widetilde{x}%
^{\prime })=g^{\prime }(\widetilde{x},\widetilde{x}^{\prime })+t.$

\begin{description}
\item[Definitions. ]  For $\widetilde{f},\widetilde{f}^{*}\in {\rm %
Diff_{t,k}^1}({\rm I}\negthinspace {\rm R^2),}$

\begin{enumerate}
\item  we say that $\widetilde{f}\leq \widetilde{f}^{*}$ if $g^{*}\leq g$
and $g^{\prime }\leq g^{*\prime }$ everywhere, where $(g,g^{\prime })$ is
associated to $\widetilde{f}$ and $(g^{*},g^{*\prime })$ is associated to $
\widetilde{f}^{*},$ as in (\ref{Tgglin}). Analogously, we say $\widetilde{f}%
\ll \widetilde{f}^{*}$ if $g^{*}<g$ and $g^{\prime }<g^{*\prime }$ everywhere;

\item  for all $r\geq 1,$ given a $C^r$ one-parameter family $(\widetilde{f}%
_t)_{t\in [a,b]}$ such that for each $a\leq t\leq b,$ $\widetilde{f}_t\in {\rm %
Diff_{t,k}^r}({\rm I}\negthinspace {\rm R^2)},$ we say the family is strongly
increasing if $\widetilde{f}_t\ll \widetilde{f}_{t^{\prime }}$ $\Leftrightarrow 
$ $t<t^{\prime }.$ We also say that a $C^r$ one-parameter family $(\widehat{f%
}_t)_{t\in [a,b]},$ such that $\widehat{f}_t\in {\rm %
Diff_{t,k}^r}({\rm T^1\times I}\negthinspace {\rm R)}$ for all $t\in[a,b]$, is strongly
increasing if it has a $C^r$ lift $(\widetilde{f}_t)_{t\in [a,b]}$ which is
strongly increasing;
\end{enumerate}
\end{description}

So 
 the one-parameter family $\widetilde{f}( 
\widetilde{x},\widetilde{y}+t)+(0,t)$ is strongly increasing. 

The next
result explains why these partial orders are important.

\begin{proposition}
\label{compte} If $(s,p,q)$ is a positive (resp. negative) triplet for $
\widetilde{f}$ and if $\widetilde{f}\leq \widetilde{f}^{*}$ (resp. $
\widetilde{f}\geq \widetilde{f}^{*}$), then $(s,p,q)$ is a positive (resp.
negative) triplet for $\widetilde{f}^{*}.$
\end{proposition}

\subsubsection{Some properties of the extremes of $\rho_V(\widehat{f})$}

The next results appeared in \cite{euerg1}:

\begin{theorem}
\label{rational0} Let $\widehat{f}\in {\rm Diff_{t,k}^1}({\rm T^1\times I}%
\negthinspace {\rm R)}$ be such that $\rho _V(\widehat{f})=[\rho
_V^{-},p/q], $ with $p/q$ a rational number. Then there exists a compact set 
$\widehat{A}\subset {\rm T^1\times I\negthinspace R}$ such that
$\widehat{f}^q(\widehat{A})-(0,p)=\widehat{A}.$
\end{theorem}
{\bf Remark:} The subset $\widehat{A}$ can be chosen as a
minimal set for $\widehat{f}^q(\bullet)-(0,p)$, but it is not true that 
$\widehat{f}^q(\bullet)-(0,p)$ always has periodic
points. Something that always holds for rational extreme points of the rotation set in 
the homotopic to the identity class.
\begin{theorem}
\label{rational1} Let $\widehat{f}\in {\rm Diff_{t,k}^1}({\rm T^1\times I}%
\negthinspace {\rm R)}$ be such that $\rho _V(\widehat{f})=[\rho
_V^{-},p/q], $ with $\rho _V^{-}<p/q,$ for some rational number $p/q.$ 
Then, assuming that the torus diffeomorphism $f$ lifted by $\widehat{f}$ has no periodic
point of vertical rotation number $p/q,$ for arbitrarily small chosen
values of $t>0,$ $\widehat{f}+(0,t)$ is the lift of a torus diffeomorphism with a $nq$%
-periodic orbit (for some $n\geq 1$) of vertical rotation number $p/q.$
And moreover, for all $t<0,$ the following holds: $\rho _V^{+}(\widehat{f}+(0,t))<p/q.$
\end{theorem}

The above theorem implies that if $\widehat{f}\in {\rm Diff_{t,k}^1}({\rm %
T^1\times I}\negthinspace {\rm R)}$ is such that $\rho _V(\widehat{f})=[\rho
_V^{-},p/q],$ and $f$ does not have periodic points with vertical rotation
number $p/q,$ then $\widehat{f}$ belongs to the boundary of $(\rho
_V^{+})^{-1}(p/q).$

\begin{theorem}
\label{irrational1} Let $\widehat{f}\in {\rm Diff_{t,k}^1}({\rm T^1\times I}%
\negthinspace {\rm R)}$ be such that $\rho _V(\widehat{f})=[\rho
_V^{-},\omega ],$ where $\omega $ is irrational. 
Then, for any $t\neq 0,$ we get that $\rho _V^{+}(\widehat{f}+(0,t))\neq \omega .$
\end{theorem}

{\bf Remark: }In order to prove the above theorem, for all $\varepsilon >0$ we
found $\widehat{f}^{*}\in {\rm Diff_{t,k}^1}({\rm T^1\times I}\negthinspace 
{\rm R),}$ such that the mappings $(g^{*},g^{*\prime })$ associated to $
\widetilde{f}^{*}$ are $\varepsilon $-$C^0$-close to $(g,g^{\prime }),$ the
mappings associated to $\widetilde{f},$ where $\widetilde{f}^{*}$ and $
\widetilde{f}$ are nearby planar lifts, respectively of $\widehat{f}^{*}$
and $\widehat{f},$ and the upper vertical rotation number has grown, that is, 
$\rho _V^{+}(\widehat{f}^{*})>\omega .$ 
Thus Proposition \ref{compte} implies that 
$\rho _V^{+}(\widehat{f}^{**})>\rho _V^{+}(\widehat{f})=\omega$ for any 
$\widehat{f}^{**}\in {\rm Diff_{t,k}^1}({\rm T^1\times I}\negthinspace {\rm R)}$
such that $\widehat{f}\ll \widehat{f}^{**}.$
 
\subsection{Prime ends compactification of open disks}

\label{sec:primends}

In this subsection we present an informal discussion on prime ends, a
subject that only appears in the proof of Theorem \ref{gerall}.

Assume $D$ is an open topological disk of an oriented surface whose boundary $\partial D$ is
not reduced to a point. 

In case $\partial D$ is a Jordan curve and $f$ is
an orientation preserving homeomorphism of that surface which satisfies $%
f(D)=D,$ it is immediate to see that $f:\partial D\rightarrow \partial D$ is
conjugate to a homeomorphism of the circle, and so a real number $\rho
(D)=rotation$ $number$ $of$ $f\mid _{\partial D}$ can be associated to this
map (up to adding an integer). Recall that, if $\rho (D)$ is rational, then
there exists a periodic point in $\partial D$ and if it is not, then there
are no such points. This is known since Poincar\'e. The difficulties arise
when we do not assume $\partial D$ to be a Jordan curve.

The prime ends compactification is a way to attach to $D$ a circle called
the circle of prime ends of $D,$ obtaining a space $D\sqcup {\rm T^1}$ with a
topology that makes it homeomorphic to the closed unit disk. If, as above,
we assume the existence of an orientation preserving homeomorphism $f$ such
that $f(D)=D,$ then $f\mid _D$ extends to $D\sqcup {\rm T^1}.$ The prime ends
rotation number of $f$ in $D,$ still denoted $\rho (D),$ is the usual
rotation number (which as before, only exists up to adding integers) of the
orientation preserving homeomorphism induced on ${\rm T^1}$ by the extension of 
$f\mid _D.$ But things may be quite different in this setting. In full
generality, it is not true that when $\rho (D)$ is rational, there are
periodic points in $\partial D$ and for some examples, $\rho (D)$ is
irrational and $\partial D$ is not periodic point free. Here we refer to 
\cite{mather} and \cite{koropatmey} for definitions, as well as to some
important theorems.

\subsection{On the existence of saddles with a full mesh}

In this subsection we present the main result of \cite{malha}. 

\begin{description}
\item[Definitions.]  If $f:{\rm T^2\rightarrow T^2}$ is a diffeomorphism,
either homotopic to the identity or to a Dehn twist, and $z\in {\rm T^2}$ is
a periodic hyperbolic saddle, we say that $z$ has a full mesh, if for any $
\widetilde{z}\in p^{-1}(z),$ $W^u(\widetilde{z})$ has a topologically
transverse intersection with $W^s(\widetilde{z})+(a,b)$ for all pairs of
integers $(a,b).$ Here $W^u(\widetilde{z})$ and $W^s(\widetilde{z})$ are the
connected components of $p^{-1}(W^u(z))$ and $p^{-1}(W^s(z))$ which contain 
$\widetilde{z}.$

We also say that $z$ has a partial mesh, if $W^u(\widetilde{z})$ has a
topologically transverse intersection with $W^s(\widetilde{z})+(a,b)$ for at
least two integer vectors $(a,b)$ which are not collinear (and thus, to
infinitely many non-collinear pairs).
\end{description}

{\bf Remark:\ }We say that a connected set $K$ has a topologically
transverse intersection with a stable manifold of a hyperbolic periodic
saddle, if there exists $z$ in this stable manifold and $r>0$ such that the
connected component of the stable manifold intersected with $B_r(z)$ which
contains $z,$ divides $B_r(z)$ into two connected components $B_{+}$ and $%
B_{-},$ and $K\cap \overline{B_r(z)}$ has a closed connected component which intersects both 
$B_{+}$ and $B_{-}$ (analogously for unstable manifolds). See for instance 
\cite{malha} and \cite{finalmente}.

\begin{theorem}
\label{malha} Let $\widehat{f}\in {\rm Diff_k^2}({\rm T^1\times I}%
\negthinspace {\rm R)\ }$(for some integer $k\neq 0)$ and suppose $p/q\in {\rm 
interior}(\rho _V(\widehat{f})).$ Then, the torus diffeomorphism $f$
lifted by $\widehat{f}$ has a periodic point of vertical rotation number $%
p/q$ which is a hyperbolic saddle with a full mesh. In the homotopic to the
identity case, if $\widetilde{f}\in {\rm Diff_0^2}({\rm I}\negthinspace {\rm %
R^2)\ }$and $(p/q,r/q)\in {\rm interior}(\rho (\widetilde{f})),$ then, the
torus diffeomorphism $f$ lifted by $\widetilde{f}$ has a hyperbolic
periodic saddle point of rotation vector $(p/q,r/q)$ with a full mesh.
\end{theorem}

If $z,w\in {\rm T^2}$ are periodic saddles, both with full or partial
meshes, then the $C^0$-$\lambda $-lemma (see for instance Proposition 1 of \cite
{london}) implies that $\overline{W^u(z)}=\overline{W^u(w)}$ and $\overline{%
W^s(z)}=\overline{W^s(w)}.$ This happens because $W^u(z)\cup W^s(z)$ and $%
W^u(w)\cup W^s(w)$ both contain closed curves in all homotopy classes. So, $%
W^u(z)$ has topologically transverse intersections with both $W^s(z)$ and $%
W^s(w),$ the same for $W^u(w)$ with respect to $W^s(z)$ and $W^s(w).$

\subsection{On a version of Pixton's Theorem to the torus}

The next result was taken from \cite{euandres}, and adapted to our notation.

\begin{theorem}
\label{comeconovo} For every integers $k$ and $r\geq 1,$ there exist
residual subsets of ${\rm Diff_k^r(T^2),}$ or in the area-preserving case of 
${\rm Diff_{k,Leb}^r(T^2),}$ denoted respectively $\chi ^r({\rm T^2})$ and $%
\chi _{{\rm Leb}}^r({\rm T^2}),$ such that: whether $f$ is homotopic to the
identity, or to a Dehn twist, assuming $f\in \chi ^r({\rm T^2})$ or $f\in
\chi _{{\rm Leb}}^r({\rm T^2})$ and $f$ has a saddle $z$ with a full mesh,
if $Y$ is a stable branch of a hyperbolic periodic saddle $p$ and $X$ is an
unstable branch of a hyperbolic periodic saddle $q$, such that $\overline{X}$
intersects $Y,$ then there exists an integer $i$ such that $Y$ intersects $%
f^i(X)$ $C^1$-transversely. Moreover, if $p$ and $q$ are in the same orbit,
then one may choose $i=0$.
\end{theorem}

{\bf Remarks:}

\begin{enumerate}
\item  {\bf \ }As explained in the Introduction after the statement of 
Theorem \ref{conjrotarea}, the subsets $\chi ^r({\rm T^2})$ and 
$\chi _{{\rm Leb}}^r({\rm T^2})$ 
consist of residual sets contained in general Kupka-Smale, or in
area-preserving Kupka-Smale $C^r$-diffeomorphisms respectively, for which
closures of stable and unstable branches of hyperbolic periodic saddles and
homoclinic and heteroclinic intersections vary continuously with
perturbations of the diffeomorphisms. See \cite{pixton}, pages 370-372.

\item  The above theorem appears in \cite{euandres} as Theorem 15. Although
an area-preserving version does not exist in that paper (because it was not
needed), as is explained in Pixton's paper \cite{pixton} it is easily
obtainable: More precisely, Theorem 8 and Lemma 9 of \cite{euandres} are
statements of results from \cite{pixton}, that are true in the
area-preserving world. And this is all we need in order to obtain
area-preserving versions of all the results in Section 3 of \cite{euandres},
culminating with Theorem \ref{comeconovo} above.
\end{enumerate}

\subsection{Generic birth and death of periodic orbits}

First we quote Theorem 20 of \cite{euandres}, which as explained in that reference, 
is essentially due to Brunovski, plus an idea of Sotomayor.

\begin{theorem}
\label{brunovski} For any $r\in \{2,3,...,\infty \},$ if $(f_t)_{t\in I}$ is
a $C^r$-generic one-parameter family of diffeomorphisms of a closed
Riemmanian manifold, then periodic points are born from only two different
types of bifurcations: saddle-nodes and period doubling. In case of
saddle-nodes, if the parameter changes, then the saddle-node unfolds into a
saddle and a sink, or source, in one direction, and the periodic point
disappears in the other direction. Moreover, at each fixed parameter, only
one saddle-node can exist.
\end{theorem}

The next result is an area-preserving version of the above one, originally
proved by Meyer \cite{meyer}, with the exception of the local picture in a
neighborhood of the saddle-elliptic point.

\begin{theorem}
\label{meyer} 
If $(f_t)_{t\in I}$ is a $C^\infty $-generic one-parameter family of
area-preserving diffeomorphisms of a closed oriented surface, then periodic
points are born from only two different types of bifurcations:
saddle-elliptic and period doubling. In case of saddle-elliptic, if the
parameter changes, then the point unfolds into a saddle and an elliptic
point (one whose eigenvalues belong to the unit circle and are not real) in
one direction and the periodic point disappears in the other direction. The
dynamics in a neighborhood of the saddle-elliptic point is as in figure 1.
\end{theorem}

\begin{figure}[!h]
	\centering
	\includegraphics[scale=0.28]{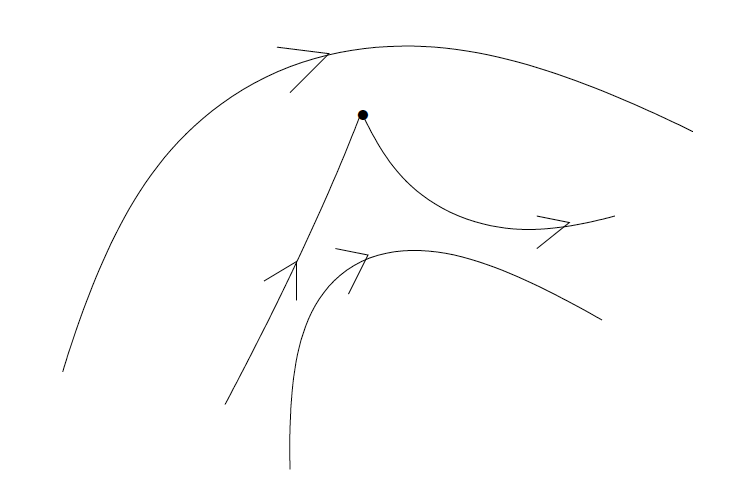}
	\caption{Diagram showing the dynamics near a saddle-elliptic periodic point in the $C^\infty$-generic case.}
\end{figure}

{\bf Remark: }
The part of the statement on the local dynamics near the saddle-elliptic
bifurcations follows from two things:

\begin{enumerate}
\item  for $C^\infty $-generic one-parameter families of area-preserving
surface diffeomorphisms, each periodic point is isolated among periodic
points of the same period, and all of them satisfy a Lojaziewicz condition,
(see \cite{lecture}, page 3 of the Summary);

\item  if $f$ is a $C^\infty $ area-preserving surface diffeomorphism with
an isolated fixed point $p$ (isolated among the fixed point set of $f$),
satisfying a Lojaziewicz condition, such that the topological index of $f$
at $p$ is zero, then both eigenvalues at $p$ are equal to $1$ and the local
dynamics is as in figure 1. A proof of this can be found in Section 2 of 
\cite{calvezeu}.
\end{enumerate}

\section{Proofs}

\subsection{Proof of Theorem 1 and Corollary 1}

{\it Sketch of the proofs.}

The main ideas used here are the following: either for every $\varepsilon >0,$ $
\widehat{f}^q(\bullet )-(0,p)$ has $\varepsilon $-pseudo orbits that move
vertically in the cylinder by arbitrarily large amounts, or not. In case for
every $\varepsilon >0$ there are such $\varepsilon $-pseudo orbits, we can find
twist maps $\widehat{h}_i\in {\rm Diff_{t,k}^r(T^1\times I\negthinspace  R)}$
arbitrarily $C^0$-close to $\widehat{f}$ such that $\rho _V^{+}(\widehat{h}%
_i)>p/q.$ And this fact, together with Lemma \ref{euerg1lemma} implies that 
$\rho _V^{+}(\widehat{f}+(0,t))>p/q$ for all $t>0.$

And in case for some $\varepsilon _0>0$ there are no such $\varepsilon _0$-pseudo
orbits, then we use the following folklore result:

\begin{proposition}
\label{folklore} Let $\widehat{h}:{\rm T^1\times I\negthinspace  %
R\rightarrow T^1\times I\negthinspace  R}$ be an orientation and
end-preserving homeomorphism. Assume that for some $a<b$ and $\varepsilon _0>0,$
there is no $\varepsilon _0$-pseudo orbit starting at some point in ${\rm %
T^1\times ]-\infty} ,a]$ and ending at some point in ${\rm T^1\times}
[b,+\infty [.$ Then, there exists a homotopically non-trivial simple closed
curve $\gamma \subset {\rm T^1\times I\negthinspace  R}$ such that $\widehat{%
h}(\gamma )\subset \gamma^-.$
\end{proposition}

{\it Proof.}

Consider the following set: $U(\varepsilon _0)=\{\widehat{z}\in {\rm %
T^1\times I\negthinspace  R:}$ $\exists $ $\varepsilon _0$-pseudo orbit for $
\widehat{h}$ starting at some point in ${\rm T^1\times ]-\infty} ,a]$ and
ending at $\widehat{z}\}.$ The hypothesis of the proposition implies that the
positively invariant, connected open set $U(\varepsilon _0)$ satisfies $%
U(\varepsilon _0)\subset {\rm T^1\times ]-\infty },b[.$ Also, if $S>0$ is
defined as $S=\sup _{\widehat{z}\in {\rm T^1\times }\{a\}\ }{\rm dist}(\widehat{h}%
( \widehat{z}),\widehat{z}),$ then $U(\varepsilon _0)\supset {\rm T^1\times
]-\infty} , a-S+\varepsilon _0].$ 

Let us show that 
\begin{equation}
\label{maisimp*}\widehat{h}({\rm closure}(U(\varepsilon _0)))\subset U(\varepsilon
_0). 
\end{equation}
To see that this concludes the proof, consider $\widehat{\Theta },$ the
connected component of the complement of ${\rm closure}(U(\varepsilon _0))$
that contains ${\rm T^1\times} ]b,+\infty [.$ Clearly, $\widehat{\Theta }%
\subset {\rm T^1\times} ]a-S+\varepsilon _0,+\infty [$ and its (connected)
boundary, $\partial \widehat{\Theta },$ separates the ends of the annulus
and is disjoint from its image under $\widehat{h},$ because $\partial 
\widehat{\Theta }\subset \partial U(\varepsilon _0)$ and $\widehat{h}(\partial
(U(\varepsilon _0)))\subset U(\varepsilon _0)$ (see expression (\ref{maisimp*})).
So, any homotopically non-trivial simple closed curve $\gamma $ contained in
the open annulus between $\partial \widehat{\Theta }$ and $\widehat{h}%
(\partial \widehat{\Theta })$ is free under $\widehat{h}$ and mapped into $\gamma^{-}.$

Thus, we are left to show that expression (\ref{maisimp*}) holds. Consider some 
$\widehat{w}\in {\rm closure} (U(\varepsilon _0)).$ Fix $\delta >0$ such that $
\widehat{h}(B_\delta (\widehat{w}))\subset B_{\varepsilon _0}(\widehat{h}( 
\widehat{w})),$ and pick a point $\widehat{w}^{*}\in B_\delta (\widehat{w})\cap
U(\varepsilon _0).$ There exists an $\varepsilon _0$-pseudo orbit 
$\{\widehat{z}_0,\widehat{z}_1,...,\widehat{z}_{n-1},\widehat{w}^{*}\}$ such
that $\widehat{z}_0\in {\rm T^1\times} ]-\infty ,a].$ As $\widehat{h}( 
\widehat{w}^{*})\in B_{\varepsilon _0}(\widehat{h}(\widehat{w})),$ the set
$\{\widehat{z}_0,\widehat{z}_1,...,\widehat{z}_{n-1},\widehat{w}^{*},\widehat{%
h}(\widehat{w})\}$ is also an $\varepsilon _0$-pseudo orbit. So, $\widehat{h}( 
\widehat{w})\in U(\varepsilon _0),$ something that implies that 
$\widehat{h}({\rm closure}%
(U(\varepsilon _0)))\subset U(\varepsilon _0).$ $\Box $

\vskip0.2truecm

These ideas when put together properly, prove Theorem 1 and Corollary 1.

\vskip0.2truecm

Back to the precise proof, as in the statements, let $k\neq 0$ be an
integer which, without loss of generality, we assume to be positive, $
\widehat{f}\in {\rm Diff_{t,k}^r(T^1\times I\negthinspace  R)}$, for any $r\geq
1,$ and let $\widetilde{f}\in {\rm Diff_{t,k}^r(I\negthinspace R^2)}$ be a lift
of $\widehat{f}$ to the plane. Clearly, 
\begin{equation}
\label{caradaf}\widetilde{f}(\widetilde{x},\widetilde{y})=(\widetilde{x}+k 
\widetilde{y}+\phi _1(\widetilde{x},\widetilde{y}),\widetilde{y}+\phi _2( 
\widetilde{x},\widetilde{y})), 
\end{equation}
where $\phi _i$ is a $1$-periodic function of $\widetilde{x}$ and $
\widetilde{y},$ for $i=1,2.$

Let us define the following constants, $k_{tw}>0,$ $A>0$ and $B>0$ as
follows: 
\begin{equation}
\label{estim1}k+\frac{\partial \phi _1}{\partial \widetilde{y}}(\widetilde{x}%
,\widetilde{y})>k_{tw}\text{ for all }(\widetilde{x},\widetilde{y})\in {\rm I%
}\negthinspace 
{\rm R^2}\text{ (twist condition)} 
\end{equation}

\begin{equation}
\label{estim2}\left| \frac{\partial \phi _1}{\partial \widetilde{x}}( 
\widetilde{x},\widetilde{y})\right| <B\text{ and }\left| \phi _2(\widetilde{x%
},\widetilde{y})\right| <A,\text{ for all }(\widetilde{x},\widetilde{y})\in 
{\rm I}\negthinspace {\rm R^2}\text{ } 
\end{equation}

The first assertion in the statement of the theorem implies that
 there exists a neighborhood of $\widehat{f}\in {\rm Diff_k^r(T^1\times I%
\negthinspace  R)}$ such that for any $\widehat{h}$ in this neighborhood, we
have $\rho _V^{+}(\widehat{h})\leq \rho _V^{+}(\widehat{f}).$ So Theorem 
\ref{irrational1} implies that $\rho _V^{+}(\widehat{f})=p/q,$ for some
rational number $p/q.$

Let us prove the following lemma:

\begin{lemma}
\label{refpediu} Under the previous hypotheses, there exists $\varepsilon _0>0$
such that there is no $\varepsilon _0$-pseudo orbit for $\widehat{f}^q(\bullet
)-(0,p):{\rm T^1\times I\negthinspace  R\rightarrow T^1\times I%
\negthinspace  R}$ starting at a point below ${\rm T^1}\times \{0\}$ and ending at
a point above ${\rm T^1}\times \{10+A+\left| p\right| +(10+B)/k_{tw}\}.$
\end{lemma}

{\it Proof.}

By contradiction, assume that for all $\varepsilon >0,$ there exist $\varepsilon $%
-pseudo orbits for $\widehat{f}^q(\bullet )-(0,p),$ starting at a point
below ${\rm T^1}\times \{0\}$ and ending at a point above ${\rm T^1}\times \{10+A+\left|
p\right| +(10+B)/k_{tw}\}.$ Denote such an $\varepsilon $-pseudo orbit by $\{ 
\widehat{z}_0,\widehat{z}_1,...,\widehat{z}_n\},$ for $n\geq 1.$ From its
choice, $p_2(\widehat{z}_0)<0$ and $p_2(\widehat{z}_n)>10+A+\left| p\right|
+(10+B)/k_{tw}.$

To the above $\varepsilon $-pseudo orbit for $\widehat{f}^q(\bullet )-(0,p),$
there corresponds the following $\varepsilon $-pseudo orbit for $\widehat{f}:$ 
\begin{equation}
\label{chainfchap}\{\widehat{z}_0,\widehat{f}(\widehat{z}_0),...,\widehat{f}%
^{q-1}(\widehat{z}_0),z_1+(0,p),...,\widehat{f}^{q-1}(\widehat{z}%
_1)+(0,p),...,
\widehat{z}_n+(0,np)\} 
\end{equation}

Modifying the points in (\ref{chainfchap}) in an arbitrarily small way, we
can obtain another $\varepsilon $-pseudo orbit for $\widehat{f}$ 
$$
\{\widehat{w}_0=(\widehat{x}_0,\widehat{y}_0),\widehat{w}_1=(\widehat{x}_1, 
\widehat{y}_1),...,\widehat{w}_{qn}=(\widehat{x}_{qn},\widehat{y}_{qn})\}, 
$$
where the $\widehat{w}_i=(\widehat{x}_i,\widehat{y}_i)$ satisfy $\widehat{x}%
_i\neq \widehat{x}_j,$ for $i\neq j.$ Clearly $p_2(\widehat{w}_0)<0$ and $%
p_2(\widehat{w}_{qn})>p.n+10+A+\left| p\right| +(10+B)/k_{tw}.$

From expression (\ref{caradaf}), 
we get that 
$$
D\widetilde{f}\mid _{(\widetilde{x},\widetilde{y})}=\left( 
\begin{array}{cc}
1+\frac{\partial \phi _1}{\partial \widetilde{x}}(\widetilde{x},\widetilde{y}%
) & k+ 
\frac{\partial \phi _1}{\partial \widetilde{y}}(\widetilde{x},\widetilde{y})
\\ \frac{\partial \phi _2}{\partial \widetilde{x}}(\widetilde{x},\widetilde{y%
}) & 1+\frac{\partial \phi _2}{\partial \widetilde{y}}(\widetilde{x}, 
\widetilde{y}) 
\end{array}
\right) 
$$
has uniformly bounded norm 
for all $(\widetilde{x},\widetilde{y})\in {\rm I}\negthinspace {\rm R^2,}$
the same holding for $\left\| D\widetilde{f}^{-1}\mid _{(\widetilde{x}, 
\widetilde{y})}\right\| .$ In ${\rm T^1\times I\negthinspace  R,}$ denote
this uniform bound for $\left\| D\widehat{f}\mid _{(\widehat{x},\widehat{y}%
)}\right\| $ and $\left\| D\widehat{f}^{-1}\mid _{(\widehat{x},\widehat{y}%
)}\right\| $ by $M>0.$ So, for all $(\widehat{x}_1,\widehat{y}_1),(\widehat{x%
}_2,\widehat{y}_2)\in {\rm T^1\times I\negthinspace  R,}$ 
$$
{\rm dist} (\widehat{f}^{\pm 1}(\widehat{x}_1,\widehat{y}_1),\widehat{f}^{\pm
1}( \widehat{x}_2,\widehat{y}_2)) <M.{\rm dist} ((\widehat{x}_1,\widehat{y}_1),
(\widehat{x}_2,\widehat{y}_2)). 
$$

Our aim now is to find some twist map $\widehat{f}^{*}\in {\rm %
Diff_{t,k}^r(T^1\times I\negthinspace  R),}$ and some of its lifts, $
\widetilde{f}^{*}\in {\rm Diff_{t,k}^r(I\negthinspace  R^2)}$ for which 
\begin{equation}
\label{ftilast}\widetilde{f}^{*}(\widetilde{x},\widetilde{y})=(\widetilde{x}%
+k\widetilde{y}+\phi _1(\widetilde{x},\widetilde{y}),\widetilde{y}+\phi
_2^{*}(\widetilde{x},\widetilde{y})), 
\end{equation}
where $\phi _2^{*}$ is $((M+1)/k_{tw}+1)\varepsilon $-$C^0$-close to $\phi _2,$ 
and for some point $\widehat{z}\in {\rm T^1\times I\negthinspace  R,}$ 
$$
p_2\left( (\widehat{f}^{*})^{nq}(\widehat{z})\right) -p_2(\widehat{z}%
)>p.n+9+A+\left| p\right| +(10+B)/k_{tw}. 
$$

As $\max _{(\widetilde{x},\widetilde{y})\in {\rm I}\negthinspace {\rm R^2}%
}\left| \phi _2^{*}(\widetilde{x},\widetilde{y})\right| <A+1$ for all $%
0<\varepsilon <((M+1)/k_{tw}+1)^{-1},$ we get from Lemma 8 of \cite{non1} that
there exists a point $\widehat{z}^{\prime }\in {\rm T^1\times I%
\negthinspace  R}$ such that one of the following possibilities holds: 
$$
\begin{array}{c}
p_2\left( ( 
\widehat{f}^{*})^{nq}(\widehat{z}^{\prime })\right) -p_2(\widehat{z}^{\prime
})>p.n+2+\left| p\right| \\ p_1\left( (\widehat{f}^{*})^{nq}(\widehat{z}%
^{\prime })\right) =p_1(\widehat{z}^{\prime }) 
\end{array}
$$
\begin{center}
or
\end{center}
$$
\begin{array}{c}
p_2\left( ( 
\widehat{f}^{*})^{nq+1}(\widehat{z}^{\prime })\right) -p_2(\widehat{z}%
^{\prime })>p.n+2+\left| p\right| \\ p_1\left( (\widehat{f}^{*})^{nq+1}( 
\widehat{z}^{\prime })\right) =p_1(\widehat{z}^{\prime }) 
\end{array}
$$
So, either the torus map $f^{*}$ has periodic points of vertical rotation
number equal to $\left( p.n+2+\left| p\right| \right) /\left( nq\right) $ or
to $\left( p.n+2+\left| p\right| \right) /\left( nq+1\right) ,$ both numbers
larger than $p/q,$ or for some integer $s,$ the triplet $(s,np+2+\left|
p\right| ,nq),$ or the triplet $(s,np+2+\left| p\right| ,nq+1)$ is positive
for $\widetilde{f}^{*}.$

But this last possibility implies that $p/q<\rho _V^{-}(\widehat{f}^{*})\leq
\rho _V^{+}(\widehat{f}^{*}).$ So, it is always the case that $p/q<\rho
_V^{+}(\widehat{f}^{*}).$

As the functions $g^{*}(\widetilde{x},\widetilde{x}^{\prime })$ and $%
g^{*\prime }(\widetilde{x},\widetilde{x}^{\prime })$ associated to $
\widetilde{f}^{*}$ will be shown to satisfy $g^{*}(\widetilde{x},\widetilde{x%
}^{\prime })=g(\widetilde{x},\widetilde{x}^{\prime })$ and 
\begin{equation}
\label{gglin*} 
\left| g^{*\prime }(\widetilde{x},\widetilde{x}^{\prime })-g^{\prime }( 
\widetilde{x},\widetilde{x}^{\prime })\right| <C^{\prime }\varepsilon ,\text{
for }C^{\prime }=(M+1)/k_{tw}+1>0, 
\end{equation}
we get from Proposition \ref{compte} and the comments right before it, that 
$$
\rho _V^{+}(\widehat{f}+(0,2C^{\prime }\varepsilon ))>p/q. 
$$

As we are assuming that $\rho _V^{+}(\widehat{h})\leq \rho _V^{+}(\widehat{f})=p/q$ 
for all $\widehat{h}\in {\rm Diff_{t,k}^r(T^1%
\times I\negthinspace  R)}$ sufficiently $C^r$-close to $\widehat{f},$  we arrived at a
contradiction because $\varepsilon>0$ could be arbitrarily small. 

So, in order to conclude the proof of the present lemma, we are left to show
the existence of $f^{*},\widehat{f}^{*},\widetilde{f}^{*}$ as above.

Recall that we are assuming the existence of an $\varepsilon $-pseudo orbit
for $\widehat{f},$ denoted 
$$
\{\widehat{w}_0=(\widehat{x}_0,\widehat{y}_0),\widehat{w}_1=(\widehat{x}_1, 
\widehat{y}_1),...,\widehat{w}_{qn}=(\widehat{x}_{qn},\widehat{y}_{qn})\}, 
$$
where $\widehat{x}_i\neq \widehat{x}_j,$ for $i\neq j$ and $p_2(\widehat{w}%
_0)<0,$ \ $p_2(\widehat{w}_{qn})>pn+10+A+\left| p\right| +(10+B)/k_{tw}.$

Clearly from the twist condition, in the vertical segment $\{\widehat{x}%
_{qn-1}\}\times [\widehat{y}_{qn-1}-\varepsilon /k_{tw},\widehat{y}%
_{qn-1}+\varepsilon /k_{tw}]$ there exists a point $\widehat{w}_{qn-1}^{\prime
} $ such that 
$$
p_1\circ \widehat{f}(\widehat{w}_{qn-1}^{\prime })=\widehat{x}_{qn}\text{
and }{\rm dist}(\widehat{f}(\widehat{w}_{qn-1}^{\prime }),\widehat{w}%
_{qn})<(M/k_{tw}+1)\varepsilon 
$$

Analogously, in the vertical segment $\{\widehat{x}_{qn-2}\}\times [\widehat{%
y}_{qn-2}-\varepsilon /k_{tw},\widehat{y}_{qn-2}+\varepsilon /k_{tw}]$ there
exists a point $\widehat{w}_{qn-2}^{\prime }$ such that 
$$
p_1\circ \widehat{f}(\widehat{w}_{qn-2}^{\prime })=\widehat{x}_{qn-1}\text{
and }{\rm dist}(\widehat{f}(\widehat{w}_{qn-2}^{\prime }),\widehat{w}%
_{qn-1}^{\prime })<((M+1)/k_{tw}+1)\varepsilon 
$$
And for $i=qn-3$ down to $0$, the situation is analogous to $i=qn-2$:
In the vertical segment $\{\widehat{x}_{i}\}\times [\widehat{%
y}_{i}-\varepsilon /k_{tw},\widehat{y}_{i}+\varepsilon /k_{tw}]$ there
exists a point $\widehat{w}_{i}^{\prime }$ such that 
$$
p_1\circ \widehat{f}(\widehat{w}_{i}^{\prime })=\widehat{x}_{i+1}\text{
and }{\rm dist}(\widehat{f}(\widehat{w}_{i}^{\prime }),\widehat{w}%
_{i+1}^{\prime })<((M+1)/k_{tw}+1)\varepsilon 
$$
So, we found a sequence of points $\widehat{w}_0^{\prime }=(\widehat{x}_0, 
\widehat{y}_0^{\prime }),\widehat{w}_1^{\prime }=(\widehat{x}_1,\widehat{y}%
_1^{\prime }),...,\widehat{w}_{qn-1}^{\prime }=(\widehat{x}_{qn-1},\widehat{y%
}_{qn-1}^{\prime }),$ $\widehat{w}_{qn}^{\prime }=\widehat{w}_{qn}=(\widehat{%
x}_{qn},\widehat{y}_{qn})$ such that, for all $i=0,1,...,qn-1:$

\begin{itemize}
\item  $p_1\circ \widehat{f}(\widehat{w}_i^{\prime })=\widehat{x}_{i+1};$

\item  $\left| \widehat{y}_i^{\prime }-\widehat{y}_i\right| <\varepsilon
/k_{tw};$

\item  ${\rm dist}(\widehat{f}(\widehat{w}_i^{\prime }),\widehat{w}_{i+1}^{\prime
})=\left| p_2(\widehat{w}_{i+1}^{\prime })-p_2(\widehat{f}(\widehat{w}%
_i^{\prime }))\right| <((M+1)/k_{tw}+1)\varepsilon ;$
\end{itemize}

With the above modification, we constructed a $((M+1)/k_{tw}+1)\varepsilon $%
-pseudo orbit $\widehat{w}_0^{\prime },\widehat{w}_1^{\prime },...,\widehat{w%
}_{qn-1}^{\prime },\widehat{w}_{qn}^{\prime }=\widehat{w}_{qn}$ such that
the image of each point is contained in the same vertical containing the
next point. In other words, in order to turn this pseudo orbit into a real
orbit for a $C^0$-nearby homeomorphism, we just have to compose $\widehat{f}$
with a diffeomorphism $\widehat{T}:{\rm T^1\times I\negthinspace  %
R\rightarrow T^1\times I\negthinspace  R}$ of the following form: 
$$
\widehat{T}(\widehat{x},\widehat{y})=(\widehat{x},\widehat{y}+\psi (\widehat{%
x})), 
$$
where $\psi :{\rm T^1 \rightarrow I\negthinspace  R}$
is a $C^\infty $ function satisfying:

\begin{itemize}
\item  $\psi (\widehat{x}_i)=p_2(\widehat{w}_i^{\prime })-p_2(\widehat{f}( 
\widehat{w}_{i-1}^{\prime }))$ for $i=1,2,...,qn;$

\item  $\left\| \psi \right\| _0<((M+1)/k_{tw}+1)\varepsilon $
\end{itemize}

Denote by $T\in {\rm Diff_0^\infty (T^2)}$ the torus diffeomorphism induced 
by $\widehat{T}.$ Also, note that any lift of $\widehat{T}$ to the plane
belongs to ${\rm Diff_0^\infty (I\negthinspace R^2).}$

If we define $f^{*}=T\circ f,$ $\widehat{f}^{*}=\widehat{T}\circ \widehat{f}$
and $\widetilde{f}^{*}=\widetilde{T}\circ \widetilde{f},$ where $\widetilde{T%
}(\widetilde{x},\widetilde{y})=(\widetilde{x},\widetilde{y}+\widetilde{\psi }%
(\widetilde{x}))$ is a lift of $\widehat{T}$ to the plane, then 
$$
\widetilde{f}^{*}(\widetilde{x},\widetilde{y})=(\widetilde{x}+k\widetilde{y}%
+\phi _1(\widetilde{x},\widetilde{y}),\widetilde{y}+\phi _2(\widetilde{x}, 
\widetilde{y})+\widetilde{\psi }(\widetilde{x}+k\widetilde{y}+\phi _1( 
\widetilde{x},\widetilde{y}))). 
$$
In the above expression, $\widetilde{\psi }:{\rm I\negthinspace  %
R\rightarrow I\negthinspace  R}$ is the $1$-periodic $C^\infty $ function
that lifts $\psi .$ So comparing the above to expression (\ref{ftilast}), we
get that $\left\| \phi _2^{*}-\phi _2\right\| _0=$ $\left\| \widetilde{\psi }%
\right\| _0<((M+1)/k_{tw}+1)\varepsilon ,$ as was stated right after expression
(\ref{ftilast}).

Moreover, 
$$
g^{*}(\widetilde{x},\widetilde{x}^{\prime })=g(\widetilde{x},\widetilde{x}%
^{\prime })\text{ and }g^{*\prime }(\widetilde{x},\widetilde{x}^{\prime
})=p_2\circ \widetilde{f}(\widetilde{x},g(\widetilde{x},\widetilde{x}%
^{\prime }))+\widetilde{\psi }(\widetilde{x}^{\prime })=g^{\prime }( 
\widetilde{x},\widetilde{x}^{\prime })+\widetilde{\psi }(\widetilde{x}%
^{\prime }). 
$$

This implies that $C^{\prime }>0$ can be expressed as a function of $%
\varepsilon >0$ exactly as in expression (\ref{gglin*}).

If $\varepsilon >0$ is small enough so that $((M+1)/k_{tw}+1)\varepsilon <1,$ then 
$$
p_2\left( (\widehat{f}^{*})^{nq}(\widehat{w}_0^{\prime })\right) -p_2( 
\widehat{w}_0^{\prime })=\widehat{y}_{qn}-\widehat{y}_0^{\prime
}>p.n+9+A+\left| p\right| +(10+B)/k_{tw}. 
$$
So, the point $\widehat{z}$ that appears right after expression (\ref
{ftilast}) is $\widehat{w}_0^{\prime }.$ And this concludes the proof of the
lemma. $\Box $

\vskip0.2truecm

The above lemma implies that assuming the first assertion in the hypotheses of the theorem, there exists $%
\varepsilon _0>0$ such that, defining 
$$
U(\varepsilon _0)=\{z\in {\rm T^1\times I\negthinspace  R:}\text{ }\exists 
\text{ }\varepsilon _0\text{-pseudo orbit for }\widehat{f}^q(\bullet )-(0,p)\text{
starting at a point in } 
$$
$$
{\rm T^1}\times ]-\infty ,0]\text{ and ending at }z\}, 
$$
then $U(\varepsilon _0)$ is open, connected, contains ${\rm T^1}\times ]-\infty
,-qA-p]$ and does not intersect ${\rm T^1}\times [10+A+\left| p\right|
+(10+B)/k_{tw},+\infty [.$ 
So Proposition \ref{folklore} implies that there exists a $\widehat{f}^q(\bullet
)-(0,p)$ free homotopically non-trivial simple closed curve $\widehat{\gamma 
}\subset {\rm T^1\times I\negthinspace  R}$ such that $\widehat{f}^q( 
\widehat{\gamma })-(0,p)\subset \widehat{\gamma }^{-}.$ 

This proves that the first assertion in the statement of the theorem implies the second. For the other 
implication, note that the second assertion implies that for any continuous map
$\widehat{f}_*$ in a $C^0$-neighborhood of $\widehat{f}$, 
$\widehat{f}_*^q(\widehat{\gamma})-(0,p)\subset \widehat{\gamma}_{-},$ and so in a $C^0$-neighborhood of
$\widehat{f},$ $\rho _V^{+}$ is smaller or equal than $p/q.$ 

In order to prove Corollary \ref{remark}, assume that $\rho _V^{-}(\widetilde{f}%
)<\rho _V^{+}(\widetilde{f})=p/q.$ 

If $\widehat{f}^q(\bullet )-(0,p)$
satisfies the (C.I.P.), then Proposition \ref{folklore} and the proof of Lemma \ref
{refpediu} imply that for all $t>0,$ $\rho _V^{+}(\widehat{f}+(0,t))>p/q,$ 
which implies that $\rho _V^{+}$ does not have a local maximum at $\widehat{f}.$

And if $\rho _V^{+}$ does not have a local maximum at $\widehat{f},$ then $
\widehat{f}^q(\bullet )-(0,p)$ does not have free homotopically non-trivial
simple closed curves, so $\widehat{f}^q(\bullet )-(0,p)$ satisfies the
(C.I.P.). This proves the equivalence in the statement of the corollary.

Finally, as $\widehat{f}^q(\bullet )-(0,p)$ is the lift to the annulus of a
torus homeomorphism homotopic to a Dehn twist, if it satisfies the C.I.P.,
then Theorem 3 of \cite{non04} implies that it has periodic orbits of all
rational rotation numbers in the annulus (with respect to the lift to the
plane $\widetilde{f}^q(\bullet )-(0,p)$). $\Box $

\vskip0.2truecm

\subsection{Proof of Theorems 2 and 3}

Let $k\neq 0$ and $r\geq 1$ be integers. The sets $O{\rm _{t,k}^r}({\rm %
T^1\times I\negthinspace  R})$ and $O{\rm _{t,k,Leb}^r}({\rm T^1\times I%
\negthinspace  R})$ that appear in the statements of Theorems \ref{conjrot} and
 \ref{conjrotarea}, defined as the subsets of maps $\widehat{f}\in {\rm %
Diff_{t,k}^r(T^1\times I\negthinspace  R)}$ or 
$\widehat{f}\in {\rm Diff_{t,k,Leb}^r(T^1\times I\negthinspace  R)}$ respectively,
such that $\rho _V$ is constant in a neighborhood of $\widehat{f},$ are both 
open by definition. We need to show that they are also dense.

If we define $O{\rm _{t,k,+}^r}({\rm T^1\times I\negthinspace  R})$ as the
subset of maps $\widehat{h}\in {\rm Diff_{t,k}^r(T^1\times I\negthinspace  R)%
}$ such that $\rho _V^{+}$ is constant in a neighborhood of $
\widehat{h},$ analogously for $O{\rm _{t,k,-}^r}({\rm T^1\times I%
\negthinspace  R})$ with respect to $\rho _V^{-},$ as both sets are clearly
open, we have:

\begin{description}
\item[Step 1.]  {\it It is enough to show that $O{\rm _{t,k,+}^r}({\rm T^1\times I%
\negthinspace  R})$ and $O{\rm _{t,k,-}^r}({\rm T^1\times I\negthinspace  R}%
) $ are both dense (the same happening in the area-preserving case).}

{\it Proof.}

This step holds trivially, because%
$$
O{\rm _{t,k}^r}({\rm T^1\times I\negthinspace  R})=O{\rm _{t,k,+}^r}({\rm %
T^1\times I\negthinspace  R})\cap O{\rm _{t,k,-}^r}({\rm T^1\times I%
\negthinspace  R}). \text{ } \Box 
$$
\end{description}

So, in order to conclude the proof, it remains to show that $O{\rm _{t,k,+}^r%
}({\rm T^1\times I\negthinspace  R})$ and $O{\rm _{t,k,Leb,+}^r}({\rm %
T^1\times I\negthinspace  R})$ are both dense (the proofs in case of 
$\rho _V^{-}$ 
are analogous). 



\begin{description}
\item[Step 2.] {\it Let $(\widehat{f}_t)_{t\in [a,b]}$ be a continuous 
one-parameter family of elements of ${\rm Diff_{t,k}^r(T^1\times I\negthinspace  
R),}$ with the additional property that it is a strongly increasing family. Then,
the map $t\mapsto \rho _V^{+}(\widehat{f}_t)$ is non-decreasing.}

{\it Proof.}

In order to prove the above statement, fix some continuous family $
\widetilde{f}_t:{\rm I\negthinspace  R^2\rightarrow I\negthinspace  R^2}$
that lifts $\widehat{f}_t,$ and note that if $t<t^{\prime },$ then $
\widetilde{f}_{t}\ll \widetilde{f}_{t^{\prime }}.$ So, if some triplet $(s,p,q)$
is non-negative for $\widetilde{f}_t,$ then it is also non-negative for $
\widetilde{f}_{t^{\prime }}.$ This easily implies that $\rho _V^{+}(\widehat{f}_{t^{\prime
}})\geq \rho _V^{+}(\widehat{f}_t).$ $\Box$

\item[Step 3.]  {\it Still under the hypothesis that $(\widehat{f}_t)_{t\in [a,b]}$ 
is a continuous one-parameter family of elements of ${\rm Diff_{t,k}^r(T^1\times I\negthinspace
R)}$ which is a strongly increasing family, if $\omega $ is irrational, 
then $t\mapsto \rho _V^{+}( 
\widehat{f}_t)$ takes the value of $\omega $ at most once.}

{\it Proof.}

This proof is contained in the proof of Theorem \ref{irrational1}
and the remark after its statement. $\Box$

%
%
%
%
%

\item[Step 4.] {\it For all $t,$ $\widehat{f}(\widehat{x},\widehat{y}+t)+(0,t)$
is a strongly increasing family.}

{\it Proof.}

See the comment right before the statement of Proposition \ref{compte}. $\Box$
\end{description}

Now let us recall that in the area-preserving case, a result by
Zehnder \cite{zehnder} implies that ${\rm Diff_{t,k,Leb}^\infty (T^2)}$ is
dense in ${\rm Diff_{t,k,Leb}^r(T^2)}$ for all $r\geq 1.$ 

And in the general case, it is an easy fact that
for all $r\geq 1$, ${\rm Diff_{t,k}^{r+1}(T^2)}$ is dense in ${\rm Diff_{t,k}^{r}(T^2)},$
moreover ${\rm Diff_{t,k}^\infty (T^2)}$ is also 
dense in ${\rm Diff_{t,k}^r(T^2)}$.   

From the above, fix some $\widehat{f}\in {\rm %
Diff_{t,k}^\infty(T^1\times I\negthinspace  R)}$ or $\widehat{f}\in {\rm %
Diff_{t,k,Leb}^\infty (T^1\times I\negthinspace  R).}$ In the remainder of
this proof, we will find a $C^\infty$-small perturbation of $\widehat{f},$ which preserves area in case 
$\widehat{f}$ preserves area, such that $\rho_V^+$ is rational and constant in a $C^0$-neighborhood 
of it, this neighborhood contained in  ${\rm Diff_{k}^0(T^1\times I\negthinspace  R)}$.
 

\begin{lemma}
\label{intover} Applying a $C^\infty$-small perturbation to 
$\widehat{f}$ if necessary, we can assume that $\rho _V(\widehat{f})$ is
either a non-degenerate interval, or a single rational number and in this
case, it is locally constant.
\end{lemma}

{\it Proof.}

In \cite{lecalvezcompte} it is proved that in case $f:{\rm T^2\rightarrow T^2%
}$ (the torus map lifted by $\widehat{f})$ does not have periodic points,
then for arbitrarily small values of $t,$ $f+(0,t)$ has periodic points, say
of vertical rotation number $r/s$ with respect to the lift $\widehat{f}%
+(0,t).$ In this case, by a $C^\infty$-small perturbation applied to $f+(0,t),$
 one of the previously obtained periodic orbits
can be made topologically non-degenerate (such a perturbation is possible
both in ${\rm Diff_{t,k}^\infty(T^2)}$ and in ${\rm Diff_{t,k,Leb}^\infty (T^2))}$%
. This implies that if we denote the perturbed mapping as $f^{\#}$ and its
lift which is close to $\widehat{f}+(0,t)$ as $\widehat{f}^{\#},$ then there
exists a $C^0$-small neighborhood of $\widehat{f}^{\#}$ in 
${\rm Diff_{k}^0(T^1 \times \rm I\negthinspace  R)}$ such that all
mappings in this neighborhood induce torus maps with periodic orbits of
vertical rotation number $r/s.$

To conclude, note that either $\rho _V$ is locally constant in a
neighborhood of $\widehat{f}^{\#}$ consisting of the single rational number $%
r/s,$ or for arbitrarily small values of $\left| s\right|,$ $\widehat{f}^{\#}+(0,s)$ has a
non-degenerate vertical rotation interval, containing $r/s$ if $\left| s\right| 
$ is sufficiently small. $\Box $

\vskip0.2truecm

So the above lemma implies that either Theorems 2 and 3 are proved (with the exception
of the part of a non-degenerate vertical rotation interval in the area-preserving case), or we
can assume that $\rho _V(\widehat{f})$ is a non-degenerate interval.

As we said before, we are left to show that both sets
$O{\rm _{t,k,+}^r}({\rm T^1\times I\negthinspace  R})$ and $O%
{\rm _{t,k,Leb,+}^r}({\rm T^1\times I\negthinspace  R})$ are dense. For
this, fix some lift $\widetilde{f}:{\rm I\negthinspace  R^2\rightarrow I%
\negthinspace  R^2}$ of $\widehat{f}.$ 
%

If $\rho _V^{+}(\widehat{f}+(0,t))$ is locally
constant for $t$ in a neighborhood of $0,$ then Lemma \ref{euerg1lemma} and 
Theorems \ref{pic} and \ref{rational1} imply that
$\rho _V^{+}(\widehat{f})$
is equal to some rational number $p/q,$ $f$ has topologically non-degenerate periodic orbits
of vertical rotation number $p/q$ and $\widehat{f}^q(\bullet )-(0,p)$ maps 
a homotopically non-trivial simple closed curve $\gamma $ into $\gamma^{-}.$ 
So the work is done: $\widehat{f}$ already belongs to 
$O{\rm _{t,k,+}^r}({\rm T^1\times I\negthinspace  R})$ or to $O%
{\rm _{t,k,Leb,+}^r}({\rm T^1\times I\negthinspace  R})$ and $\rho _V^{+}$ is 
constant and rational in a neighborhood of $\widehat{f}$ in 
${\rm Diff_{k}^0(T^1\times I\negthinspace  R)}.$

Assume now that $\rho _V^{+}(\widehat{f}+(0,t))$ is not locally constant in any
neighborhood of $0.$ Fix some $\eta >0,$ an arbitrarily small number
such $\rho _V(\widehat{f}+(0,t))$ is a non-degenerate interval for all $%
\left| t\right| <\eta $ (the $C^0$-continuity of $\rho _V$ implies that this
is possible). Either when $f$ preserves area or not, for any
$0<\eta_1<\eta$ and $\delta >0,$ we can find a $C^\infty$-generic family 
$(\widehat{f}_t^{*})_{t\in [-\eta_1,\eta_1]},$ 
lifted by $\widetilde{f}_t^{*},$ $\delta $-$C^\infty$-close
to $\widehat{f}(\widehat{x},\widehat{y}+t/2)+(0,t/2)$ and $\widetilde{f}( 
\widetilde{x},\widetilde{y}+t/2)+(0,t/2)$ respectively, for all 
$t\in [-\eta_1,\eta_1],$ satisfying some special properties, see Theorems \ref
{brunovski} and \ref{meyer}. 

\begin{description}
\item[Step 5.]  {\it If $\delta >0$ is sufficiently small, then $(\widehat{f}%
_t^{*})_{t\in [-\eta_1,\eta_1]}$ is strongly increasing.}

{\it Proof.}

For each $t\in [-\eta_1,\eta_1],$ $\widetilde{f}_t^{*}$ is $\delta $-$C^\infty$%
-close to $\widetilde{f}(\widetilde{x},\widetilde{y}+t/2)+(0,t/2).$ So the
pair of mappings $(g_t^{*},g_t^{*\prime })$ associated to $\widetilde{f}%
_t^{*}$ satisfy the following inequalities:%
$$
\left| \frac \partial {\partial t}g_t^{*}(\widetilde{x},\widetilde{x}%
^{\prime })-(-1/2)\right| <C(\delta )\text{ and }\left| \frac \partial
{\partial t}g_t^{*\prime }(\widetilde{x},\widetilde{x}^{\prime })-1/2\right|
<C(\delta ), 
$$
for all $(\widetilde{x},\widetilde{x}^{\prime })\in {\rm I\negthinspace  R^2}
$ and some constant $C(\delta )\rightarrow 0$ as $\delta \rightarrow 0.$ So,
if $\delta >0$ is small enough, $\frac \partial {\partial t}g_t^{*}( 
\widetilde{x},\widetilde{x}^{\prime })<0$ and $\frac \partial {\partial
t}g_t^{*\prime }(\widetilde{x},\widetilde{x}^{\prime })>0$ for all $( 
\widetilde{x},\widetilde{x}^{\prime })\in {\rm I\negthinspace  R^2,}$ that
is, $(\widehat{f}_t^{*})_{t\in [-\eta_1,\eta_1]}$ is a strongly increasing
family. We are using the fact that if the pair $(g,g^{\prime })$ is
associated to $\widetilde{f},$ then $(g-t/2,g^{\prime }+t/2)$ is associated
to $\widetilde{f}(\widetilde{x},\widetilde{y}+t/2)+(0,t/2).$ $\Box$
\end{description}

From our hypotheses, no matter the value of $\eta_1>0,$ 
$\rho _V^{+}(\widehat{f}-(0,\eta_1))<\rho _V^{+}(\widehat{f}+(0,\eta_1)).$ 
Fix $0<\eta_1<\eta$ and $\delta >0$ small enough (among other
requirements, $\delta <\eta_1 /4$), so that 
$\rho _V^{+}(\widehat{f}_{-\eta_1}^{*})<\rho _V^{+}(\widehat{f}_{\eta_1}^{*}),$ 
$\widehat{f}_t^{*}$ is 
$\eta$-$C^r$-close to $\widehat{f}$ for all $t\in [-\eta_1,\eta_1],$ 
$( \widehat{f}_t^{*})_{t\in [-\eta_1,\eta_1]}$ is a strongly increasing
family and $\rho _V(\widehat{f}_t^{*})$ is a non-degenerate interval for all 
$t\in [-\eta_1,\eta_1].$ 

Pick a rational $\rho _V^{+}(\widehat{f}%
_{-\eta_1}^{*})<p/q<\rho _V^{+}(\widehat{f}_{\eta_1}^{*})$ and  
$t^{\prime }\in ]-\eta_1,\eta_1[$ such that
$$
t^{\prime }=\inf \{t\in [-\eta_1,\eta_1]:\rho _V^{+}(\widehat{f}_t^{*})\geq p/q\}. 
$$
Clearly, $-\eta_1<t^{\prime }<\eta_1,$ $\rho _V^{+}(\widehat{f}%
_{t^{\prime }}^{*})=p/q,$ $\rho _V^{+}(\widehat{f}_t^{*})\geq p/q$ for all $%
t>t^{\prime },$ and for $t<t^{\prime },$ $\rho _V^{+}(\widehat{f}%
_t^{*})<p/q. $

Theorem \ref{rational1} implies that for some $t^{\prime \prime }\geq
t^{\prime },$ $\rho _V^{+}(\widehat{f}_{t^{\prime \prime }}^{*})=p/q$ and 
$f_{t^{\prime \prime }}^{*}$ has finitely many $nq$-periodic points for some
integer $n>0$ with vertical rotation number $p/q$ which are all
degenerate: saddle-node in the general case, or saddle-elliptic in the
area-preserving case (this follows from Theorems \ref{brunovski} and \ref
{meyer}). Moreover, we can assume, decreasing $t^{\prime \prime }$ if
necessary, that $f_{t^{\prime \prime }}^{*}$ has no $iq$-periodic points of
vertical rotation number $p/q$ for $i=1,2,...,n-1.$ This 
follows from the fact that $(\widehat{f}_t^{*})_{t\in [-\eta_1,\eta_1]}$ is
a strongly increasing family. 

As the $nq$-periodic points with vertical 
rotation number $p/q$ that exist for $f_{t^{\prime \prime }}$ are all degenerate,
the choice of $n$ implies that, for $t<t^{\prime \prime },$  
$f_t^{*}$ does not have periodic points of vertical rotation number $p/q$ and 
period smaller or equal to $nq,$ and 
for all $t-t^{\prime \prime }>0$ sufficiently small, $f_t^{*}$ has
hyperbolic $nq$-periodic saddles of vertical rotation number $p/q,$ as the
degenerate $nq$-periodic points for $f_{t^{\prime \prime }}^{*}$ bifurcate in the way explained in Theorems 
\ref{brunovski} and \ref{meyer}. In this way, for any $t-t^{\prime \prime }>0$
sufficiently small, there exists a $C^0$-neighborhood $V$ of 
$\widehat{f}_t^{*}$ in ${\rm Diff_{k}^0(T^1\times I\negthinspace  R)},$ such that all $
\widehat{h}\in V$ are lifts of torus maps with $nq$-periodic points
of vertical rotation
number $p/q,$ because saddle periodic points are topologically non-degenerate, 
they cannot be destroyed by $C^0$-small perturbations.
Thus the upper extreme of the vertical rotation interval cannot decrease, that is,  $\rho _V^{+}(\widehat{h})\geq p/q.$

The local dynamics, either in the saddle-node, or in the saddle-elliptic
case are what we called {\it trivializable} in Section 3 of \cite{calvezeu}.
This means that there exists an isolating neighborhood of the $nq$-periodic
point in question, such that, after a local change of coordinates, all
points in this neighborhood, with the exception of the periodic point
itself, move to the right under iterations of the dynamics which are
multiples of $nq,$ see Figure 2.

\begin{figure}[!h]
	\centering
	\includegraphics[scale=0.48]{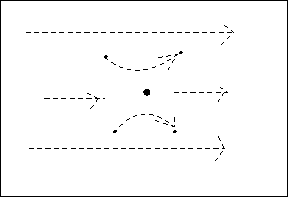}
	\caption{Diagram showing the dynamics in a neighborhood of a {\it trivializable} 
periodic point, after an appropriate coordinate change.}
\end{figure}

The proofs that saddle-node and saddle-elliptic periodic points are
{\it trivializable} appear in Propositions 6 and 8 of that paper, respectively
(although Proposition 8 asks for area-preservation and real analiticity, the
generic family that comes from Theorem \ref{meyer} satisfies everything
needed there, namely, for each integer $n>0,$ there are finitely many $n$-periodic 
points, and each
periodic point satisfies a Lojaziewicz condition).

Now we conclude our proof using Proposition 9 of \cite{calvezeu}. Adapted to our 
setting, it says
that there exists a $C^0$-neighborhood $W$ of $\widehat{f}_{t^{\prime \prime
}}^{*}$ in ${\rm Diff_{k}^0(T^1 \times \rm I\negthinspace  R)}$ such that for all 
$\widehat{h}\in W,$ $\rho _V^{+}(\widehat{h})\leq
p/q.$ So, if we pick $t>t^{\prime \prime }$ such that $\widehat{f}_t^{*}\in
W,$ we get that for all $\widehat{h}\in V\cap W,$ $\rho _V^{+}(\widehat{h}%
)=p/q.$

Although Proposition 9 of \cite{calvezeu} is stated for the homotopic to the
identity class, it also holds for maps homotopic to Dehn twists. 

Its proof
is based on a series of local deformations and perturbations: To be more precise, if by
contradiction, the adapted statement of the proposition did not hold, then there
would be 
a deformation $\widehat{h}\in {\rm Diff_k^0(T^1\times I\negthinspace  R)}$ 
of the map $\widehat{f}_{t^{\prime \prime }}^{*},$ 
such that $\rho_V^{-}(\widehat{h})<p/q<\rho_V^{+}(\widehat{h}),$ this deformation supported 
outside neighborhoods of all the $nq$-periodic points of $f_{t^{\prime \prime }}^{*}$ 
with vertical rotation number $p/q.$ Moreover, if $h$ is the torus map lifted by 
$\widehat{h}$, then the set of $h$-periodic points with vertical rotation 
number $p/q$ and period smaller or equal to $nq$ is equal to the same set for 
$f_{t^{\prime \prime }}^{*},$ in other words, the only possible period is $nq.$
 
This cannot be achieved in general, 
the proof of Proposition 9 of [7] strongly uses the fact that the dynamics is 
{\it trivializable} in neighborhoods of all the $nq$-periodic points with vertical rotation number $p/q.$  

The remaining part of the argument is to perturb $h$ inside the {\it trivializable}
neighborhoods in order to erase the 
$nq$-periodic points of vertical rotation number $p/q.$ 
Denote this perturbation of $h$ by $h^\prime$.

The above construction is possible by the choice of $n$ and $t^{\prime \prime
}$ above, and the fact that it is very easy to destroy {\it trivializable} periodic points:
any adequate arbitrarily $C^0$-small local translation supported inside the 
{\it trivializable} neighborhoods will do the job. 
So, as $\rho _V$ varies continuously in
the $C^0$-topology, $p/q$ would still be an interior point of $\rho _V(\widehat{h}^\prime),$ 
and $p/q$ would not be realized by a $q$-periodic orbit, as stated in
Theorem 5.3 of \cite{alunofranks}. This contradiction shows the existence of
the open neighborhood $W\subset {\rm Diff_k^0(T^1\times I\negthinspace  R)\ }$of $
\widehat{f}_{t^{\prime \prime }}^{*}$ such that for all $\widehat{h}\in W,$ $%
\rho _V^{+}(\widehat{h})\leq p/q,$ as explained above.

Summarizing, we found $\widehat{h}^{\#\#}\in {\rm Diff_{t,k}^\infty(T^1%
\times I\negthinspace  R)}$ or $\widehat{h}^{\#\#}\in {\rm Diff_{t,k,Leb}^\infty (T^1\times I%
\negthinspace  R),}$ such that $\widehat{h}^{\#\#}$ equals $\widehat{f}$ in case $%
\rho _V^{+}$ is locally constant in a neighborhood of $\widehat{f},$ or 
$\widehat{h}^{\#\#}$ equals $\widehat{f}_t^{*}$ if it is not, such that:

\begin{enumerate}
\item  $\widehat{h}^{\#\#},\widetilde{h}^{\#\#}$ are $\eta$-$C^r$-close
to $\widehat{f},\widetilde{f}$

\item  $\rho _V^{+}(\widehat{h}^{\#\#})=p/q$ and $\rho _V^{+}$ is locally
constant in a $C^0$-neighborhood of $\widehat{h}^{\#\#}$ in 
${\rm Diff_k^0(T^1\times I\negthinspace  R)\ }.$
\end{enumerate}

We are left to show that, in the area-preserving case, $\rho _V^{-}<\rho
_V^{+}.$ For this, just consider the vertical rotation number of the
Lebesgue measure (i.e. area).

Namely, for a fixed $\widehat{f}\in {\rm Diff_{t,k,Leb}^r(T^1\times I%
\negthinspace  R)}$ such that $\rho _V(\widehat{f})=[r/s,p/q]$ and the
vertical rotation interval is locally constant in a neighborhood of $
\widehat{f},$ note that Lebesgue measure's vertical rotation number of $
\widehat{f}+(0,s),$ denoted $\rho _{V,{\rm Leb}}(\widehat{f}+(0,s))$ is
equal to $\rho _{V,{\rm Leb}}(\widehat{f})+s$ (for all $s\in {\rm I}%
\negthinspace 
{\rm R).}$ As $\rho _V(\widehat{f}+(0,s))=\rho _V(\widehat{f})=[r/s,p/q]$
for all $\left| s\right| $ sufficiently small, and $\rho _{V,{\rm Leb}}( 
\widehat{f}+(0,s))\in \rho _V(\widehat{f}+(0,s))=[r/s,p/q],$ also for $%
\left| s\right| $ sufficiently small, it must be the case that $r/s<p/q.$
And the proof is over. $\Box $

\vskip0.2truecm

\subsection{A technical result on the existence of partial mashes}

The next result, Theorem \ref{gerall},\ref{gerall}$^\prime$, 
is the only one in this paper that also applies to the
identity homotopy class. It is auxiliary to Theorem \ref{locconst} and
Corollary \ref{locconsttran}, but it has its own interest.

In \cite{finalmente}, we proved a version of Theorem \ref{gerall} for
area-preserving homotopic to the identity maps, and as a consequence, we
showed that whenever the rotation set of a generic one-parameter family of
area-preserving diffeomorphisms of the torus homotopic to the identity
changes as the parameter changes, then certain homoclinic tangencies in the
torus (heteroclinic in the plane) appear and unfold generically, giving rise
to all phenomena associated to such an unfolding, birth of elliptic islands
(infinitely many), etc.

Below we present a simpler proof in the area-preserving case, and also one
that works in general. So Theorem \ref{gerall} implies a version of the main result of \cite
{finalmente} to any generic one-parameter family of diffeomorphisms of the
torus, either homotopic to Dehn twists or to the identity.

More precisely, whenever the rotation set (or vertical rotation interval,
in the Dehn twist case) of a generic one-parameter family of diffeomorphisms
of the torus changes as the parameter changes, then certain homoclinic
tangencies unfold generically. These tangencies are between stable and
unstable manifolds of hyperbolic periodic saddles of rational rotation
vectors (or numbers) which are eaten by the rotation set as the parameter
changes (by eaten, we mean that as the parameter changes, these points
modify their status, from boundary to interior points).

\begin{theorem}
\label{gerall} (Homotopic to Dehn twists version) For any integers $k\neq 0$ 
and $r\geq 1,$ let $f\in {\rm Diff_k^r(T^2)\cap \chi ^r(T^2)}$ or 
$f\in {\rm Diff_{k,Leb}^r(T^2)\cap \chi
_{Leb}^r(T^2)}$  be such 
that the vertical rotation interval 
$\rho _V(\widehat{f})$ has interior, for some fixed lift 
$\widehat{f}\in {\rm Diff_k^r(T^1\times I\negthinspace  R)}$ of $f$. 
Then for any rational number $\rho 
$ in the boundary of the vertical rotation interval such that $f$ has $%
\rho $-periodic orbits, there exist $\rho $-hyperbolic periodic saddles for 
$f$ with a partial mesh.
\end{theorem}

\noindent {\bf Theorem \ref{gerall}$^\prime$} {\it (Homotopic to the identity version) For any $r\geq 1,$ 
let $f\in {\rm Diff_0^r(T^2)\cap \chi ^r(T^2)}$ or 
$f\in {\rm Diff_{0,Leb}^r(T^2)\cap \chi_{Leb}^r(T^2)}$ 
be such that the rotation set $\rho (\widetilde{f})$ 
has interior, for some fixed some lift 
$\widetilde{f}\in {\rm Diff_0^r(I\negthinspace  R^2)}$ of $f$. 
Then for any rational rotation vector $\rho$ 
in the boundary of the rotation set such that $f$ has $\rho $-periodic orbits, 
there exist $\rho $-hyperbolic periodic saddles for 
$f$ with a partial mesh.}

{\bf Remarks:\ } 

\begin{itemize}

\item If the vertical rotation interval has interior for some lift of a homeomorphism 
of the torus homotopic to a Dehn twist, then it has interior for all lifts, same thing 
happening in the homotopic to the identity case.

\item Recall that $\chi^r({\rm T^2})$ and $\chi_{Leb}^r({\rm T^2})$ come from
Theorem \ref{comeconovo}. 

\item By saying that $f$ has a $\rho $-periodic orbit, we mean that
it has a periodic orbit of vertical rotation number $\rho $ (or rotation
vector $\rho ,$ in case $k=0$).

\end{itemize}

{\it Proof.}

The proof is almost identical in case $k=0$ or $k\neq 0$. We will assume that
$k\neq 0$ in the writing and stress the differences when necessary.

If $r\geq 2,$ then Theorem \ref{malha} implies that for all $p/q\in {\rm interior}(\rho
_V(\widehat{f})),$ there exists a $f$-periodic hyperbolic saddle in the
torus whose vertical rotation number is $p/q,$ which has a full mesh. If $%
r=1,$ recall that ${\rm Diff_k^2(T^2)}$ is dense in ${\rm Diff_k^1(T^2),}$
and ${\rm Diff_{k,Leb}^2(T^2)}$ is also dense in ${\rm Diff_{k,Leb}^1(T^2),}$
again see \cite{zehnder}. As we are assuming $C^r$-generic conditions, 
topologically transverse intersections are just $C^1$-transverse intersections,
which are stable under $C^1$-small perturbations. And as the vertical rotation intervals
vary continuously in the $C^0$ topology, there exist open and dense
subsets, $O_{f.m.}^1\subset \{h\in {\rm Diff_k^1(T^2)}:\rho _V(\widehat{h})$
has interior for $\widehat{h}\in {\rm Diff_k^1(T^1\times I\negthinspace  R)}$
that lifts $h\}$ and $O_{{\rm Leb},f.m.}^1\subset \{h\in {\rm %
Diff_{k,Leb}^1(T^2)}:\rho _V(\widehat{h})$ has interior for $\widehat{h}\in 
{\rm Diff_{k,Leb}^1(T^1\times I\negthinspace  R)}$ that lifts $h\}$ such that every 
$h\in O_{f.m.}^1$ or $h\in O_{{\rm Leb},f.m.}^1$ has hyperbolic periodic saddles
with full mesh. So, when $r=1,$ apart from what was explained in Subsection
2.4 on the definition of $\chi ^1({\rm T^2})$ and $\chi _{{\rm Leb}}^1({\rm %
T^2}),$ we also assume that $\chi ^1({\rm T^2})\subset O_{f.m.}^1$ and $\chi
_{{\rm Leb}}^1({\rm T^2})\subset O_{{\rm Leb},f.m.}^1.$ An analogous situation 
holds for $k=0$.

In order to prove the present result, assume that $p/q\in \partial \rho _V( 
\widehat{f})$ and the induced torus map $f$ that belongs to ${\rm %
Diff_k^r(T^2)\cap \chi ^r(T^2)}$ or ${\rm Diff_{k,Leb}^r(T^2)\cap \chi
_{Leb}^r(T^2),}$ has $q$-periodic orbits with vertical rotation
number $p/q.$ 
As the Euler characteristic
of ${\rm T^2}$ is zero, and in both the identity and the Dehn twist homotopy
classes, both eigenvalues of the action of $f$ on the first homology group
of the torus are equal to $1,$ a Lefschetz--Nielsen type theorem (see \cite
{katok}) implies that (when $k\neq 0$)
\begin{equation}
\label{lefschetz}\sum_{z\in {\rm Fix}(\widehat{f}^q(\bullet )-(0,p))}{\rm ind}(f^q\mid
_z)=0, 
\end{equation}
where ${\rm Fix}(\widehat{f}^q(\bullet )-(0,p))=\{z\in {\rm T^2:}$ $f^q(z)=z$ and $%
z$ has rot. number $p/q\},$ a similar condition holding in the homotopic to
identity case.

As $f\in \chi ^r({\rm T^2})$ or $f\in \chi _{{\rm Leb}}^r({\rm T^2}),$ 
${\rm Fix}(\widehat{f}^q(\bullet )-(0,p))$ is finite (for any rational $p/q$), and the
topological index of $f^q$ on a $q$-periodic point only assumes the values $-1$ or $1.$
As $-1$ corresponds to hyperbolic saddles with positive eigenvalues, our
hypotheses imply that ${\rm Fix}(\widehat{f}^q(\bullet )-(0,p))$ contains $n$ (for
some $n\geq 1$) hyperbolic $q$-periodic orbits of saddle type, each one of
topological index $-1$ (all of them with vertical rotation number $p/q),$ denoted 
\begin{equation}
\label{conjtotal}\{Q_1^1,Q_2^1,...,Q_q^1\},\{Q_1^2,Q_2^2,...,Q_q^2\},...,%
\{Q_1^n,Q_2^n,...,Q_q^n\}. 
\end{equation}
Each $\{Q_1^i,Q_2^i,...,Q_q^i\}$ is a single $q$-periodic orbit.

If for some $Q_j^i,$ both $W^u(Q_j^i)$ and $W^s(Q_j^i)$ are unbounded when
lifted to the plane (meaning that each connected component of their lifts to
the plane is unbounded), then:

\begin{itemize}
\item  for all $1\leq l\leq q,$ $W^u(Q_l^i)$ and $W^s(Q_l^i)$ are unbounded
when lifted to the plane;

\item  as $f$ has a saddle with a full mesh, we get that $Q_l^i$ (for all $%
1\leq l\leq q$) has a partial mesh (a proof of this when $k=0$ appears in Subsection
2.1 of \cite{finalmente} and works when $k\neq 0$ as well);
\end{itemize}

So, all we need to show, is that there exists $i\in \{1,...,n\}$ and $1\leq
j\leq q$ such that both $W^u(Q_j^i)$ and $W^s(Q_j^i)$ are unbounded when
lifted to the plane. As explained above, if this happens for some $i$ and $%
j, $ then it happens for all $Q_l^i,$ $1\leq l\leq q.$

By contradiction, suppose this is not the case. In other words, assume that
for every $i\in \{1,...,n\}$ and all $1\leq l\leq q,$ $W^u(Q_l^i)$ or $%
W^s(Q_l^i)$ is bounded when lifted to the plane. We will show that this
violates expression (\ref{lefschetz}).

Order the periodic orbits in (\ref{conjtotal}), in the following way: for
some $0\leq s^{*}\leq n,$ and for all $1\leq i\leq s^{*},$ $W^u(Q_l^i)$ (for
all $1\leq l\leq q)$ is bounded when lifted to the plane, and for all $%
s^{*}+1\leq i\leq n,W^s(Q_l^i)$ (for all $1\leq l\leq q)$ is bounded when
lifted to the plane and $W^u(Q_l^i)$ is not. In other words, $s^{*}<n$
implies that for all $s^{*}+1\leq i\leq n$ and $1\leq l\leq q,$ $W^u(Q_l^i)$
is unbounded when lifted to the plane; if for some periodic orbit in (\ref
{conjtotal}), both its stable and unstable manifolds are bounded when lifted
to the plane, then its index $i$ is smaller or equal to $s^{*}$ (clearly, if $%
s^{*}=0,$ then for all $1\leq i\leq n$ and all $1\leq l\leq q,$ $W^u(Q_l^i)$
is unbounded when lifted to the plane).

The above assumption implies that the vertical rotation number of any
point in the $f$-invariant, closed subset $\theta \subset {\rm T^2}$ given
by 
\begin{equation}
\label{defteta}\theta =\left( \cup _{i=1}^{s^{*}}\cup _{l=1}^q\overline{%
W^u(Q_l^i)}\right) \cup \left( \cup _{i=s^{*}+1}^n\cup _{l=1}^q\overline{%
W^s(Q_l^i)}\right) , 
\end{equation}
is equal to $p/q.$

As we are assuming that the vertical rotation interval has interior, we
get that the complement of ${\rm Filled}(\theta )$ is fully essential.

The proof of the present theorem in the area-preserving case is much easier,
and is concluded as follows. Let $M\subset {\rm T^2}$ be a connected
component of ${\rm Filled}(\theta ).$ It is $f^q$-invariant and $M^c$ is
connected and fully essential. The prime ends rotation number of $f^q$
restricted to the boundary of $M$ is irrational (and thus, different from
zero), see \cite{mather}. So, the orbit of a point outside, but close to $M,$
is turning around it, and thus the topological index of $f^q$ restricted to
any (contractible) simple closed curve $\alpha $ sufficiently close to $%
\partial $$M$ in the Hausdorff topology, ${\rm interior}(\alpha )\supset M,$
is equal to $1.$

This means that $\sum_{z\in Fix(\widehat{f}^q(\bullet )-(0,p))\cap
M}{\rm ind}(f^q\mid _z)=1.$ As ${\rm Filled}(\theta )$ is non-empty and it
contains all the periodic orbits in (\ref{conjtotal}), the ones with
negative topological indices, the proof in the area-preserving case is over: a
contradiction with expression (\ref{lefschetz}) was found.

From now on, assume $f\in {\rm Diff_k^r(T^2)\cap \chi ^r(T^2).}$ As in the
area-preserving case, we are going to show that 
$$
\sum_{z\in {\rm Fix}(\widehat{f}^q(\bullet )-(0,p))\cap \theta ^{\prime
}}{\rm ind}(f^q\mid _z)>0, 
$$
for some subset $\theta ^{\prime }\supset \theta ,$ and thus, containing all
periodic orbits in (\ref{conjtotal}). As these are the totality of $q$%
-periodic points of vertical rotation number $p/q$ and topological index $-1,$ we
again get a contradiction with expression (\ref{lefschetz}).

An important remark is the following: for all $1\leq i,i^{\prime }\leq n$
and $1\leq l,l^{\prime }\leq q,$ Theorem \ref{comeconovo} implies that if $%
Q_{l^{\prime }}^{i^{\prime }}\in \overline{W^u(Q_l^i)},$ then $\overline{%
W^u(Q_{l^{\prime \prime }}^{i^{\prime }})}\subset \overline{W^u(Q_l^i)},$
for some $1\leq l^{\prime \prime }\leq q,$ an analogous result holding for
the closure of stable manifolds.

\begin{lemma}
\label{gerale} In the above setting, if $\cup _{l=1}^m\overline{%
W^u(Q_{j_l}^{i_l})}$ (for some integer $m\geq 1$ and sequences $1\leq
i_l\leq s^{*},$ $1\leq j_l\leq q)$ is connected, then

$$
\sum_{z\in {\rm Fix}(\widehat{f}^q(\bullet )-(0,p))\cap {\rm Filled}(\cup _{l=1}^m 
\overline{W^u(Q_{j_l}^{i_l})})}{\rm ind}(f^q\mid _z)=1, 
$$
a similar result holding for stable manifolds whose union is connected and
whose lifts to the plane are bounded.
\end{lemma}

{\it Proof.}

From our hypotheses, ${\rm Filled}(\cup _{l=1}^m\overline{W^u(Q_{j_l}^{i_l})}%
)$ has a connected, fully essential complement. As we did in the
area-preserving case, if the prime ends rotation number of $f^q$ restricted
to the boundary of ${\rm Filled}(\cup _{l=1}^m\overline{W^u(Q_{j_l}^{i_l})})$
is non-zero, then the sum of topological indices that appear in the statement of the
lemma must be $1.$

So, in order to finish the proof, we are left to understand what happens
when the prime ends rotation number described above is zero. In this case,
Theorem D of \cite{euandres}, implies that ${\rm Filled}(\cup _{l=1}^m 
\overline{W^u(Q_{j_l}^{i_l})})$ is an attractor, that is, there exists a
contractible simple closed curve $\gamma \subset {\rm T^2,}$ ${\rm interior}%
(\gamma )\supset {\rm Filled}(\cup _{l=1}^m\overline{W^u(Q_{j_l}^{i_l})})$
such that $f^q(\gamma )\subset {\rm interior}(\gamma )$ and $\cap
_{i=0}^{+\infty }$$f^{qi}({\rm interior}(\gamma ))={\rm Filled}(\cup
_{l=1}^m \overline{W^u(Q_{j_l}^{i_l})}).$ So, 
$$
\sum_{z\in {\rm Fix}(\widehat{f}^q(\bullet )-(0,p))\cap {\rm Filled}(\cup _{l=1}^m 
\overline{W^u(Q_{j_l}^{i_l})})}{\rm ind}(f^q\mid _z)={\rm ind}(f^q\mid _\gamma )=1, 
$$
and the proof is over. $\Box $

{\bf Remark: }In Theorem D of \cite{euandres}, only the case $m=1$ is
considered, but the proof for the general case, as long as a CONNECTED set
of the form $\cup _{l=1}^m\overline{W^u(Q_{j_l}^{i_l})}$ is considered, is
exactly the same. This is a place where Theorem \ref{comeconovo} of Section
2.4 is used: it is essential in the proof of Theorem D of \cite{euandres}.

\vskip0.2truecm

The above lemma implies that the theorem is proved when $s^{*}=0$ or $s^{*}=n.$ 
In order to see this, assume $s^{*}=n.$ In this case, each connected 
component of ${\rm Filled}(\theta)$ is of the form 
${\rm Filled}(\cup _{l=1}^m \overline{W^u(Q_{j_l}^{i_l})}),$ for some 
integer $m\geq 1$ and sequences $1\leq i_l\leq n,$ $1\leq j_l\leq q$ such that
 $\cup _{l=1}^m\overline{W^u(Q_{j_l}^{i_l})}$ is connected. So, the previous lemma 
implies that the sum of topological
indices of $f^q$ at all points in ${\rm Fix}(\widehat{f}^q(\bullet )-(0,p))\cap 
{\rm Filled}(\theta)$ is positive, a contradiction as we already explained.
The case $s^{*}=0$ is analogous, one just has to consider stable manifolds instead. 

Thus, suppose $0<s^{*}<n.$

\begin{proposition}
\label{propsab0} As in the statement of the above lemma, suppose that $\cup
_{l=1}^m\overline{W^u(Q_{j_l}^{i_l})}$ (for some integer $m\geq 1$ and
sequences $1\leq i_l\leq s^{*},$ $1\leq j_l\leq q)$ is connected. Then, for
any $s^{*}<i^{\prime }\leq n$ and $1\leq l^{\prime }\leq q,$ $%
W^s(Q_{j^{\prime }}^{i^{\prime }})$ is contained in a single connected
component of $(\cup _{l=1}^m\overline{W^u(Q_{j_l}^{i_l})})^c.$
\end{proposition}

{\it Proof. }

Otherwise, for some $1\leq i\leq s^{*}$ and $1\leq j\leq q,$ $\overline{%
W^u(Q_j^i)}$ would intersect $W^s(Q_{j^{\prime }}^{i^{\prime }})$ and so,
Theorem \ref{comeconovo} would imply that $W^u(Q_j^i)$ intersects $%
W^s(Q_{j^{\prime \prime }}^{i^{\prime }}),$ for some $1\leq j^{\prime \prime
}\leq q.$ And thus, it accumulates on $W^u(Q_{j^{\prime \prime }}^{i^{\prime
}}).$ And this is a contradiction, since $W^u(Q_{j^{\prime \prime
}}^{i^{\prime }})$ is unbounded when lifted to the plane and $W^u(Q_j^i)$ is
not. In particular, $Q_{j^{\prime }}^{i^{\prime }}\notin \cup _{l=1}^m 
\overline{W^u(Q_{j_l}^{i_l})}.$ $\Box $

\vskip0.2truecm

{\bf Remark: }Similarly, if $\cup _{l=1}^{m^{\prime }}\overline{%
W^s(Q_{j_l^{\prime }}^{i_l^{\prime }})}$ (for some integer $m^{\prime }\geq
1 $ and sequences $s^{*}+1\leq i_l^{\prime }\leq n,$ $1\leq j_l^{\prime
}\leq q)$ is connected, then for any $1\leq i\leq s^{*}$ and $1\leq j\leq q,$
$W^u(Q_j^i)$ is contained in a single connected component of $(\cup _{l=1}^m 
\overline{W^s(Q_{j_l^{\prime }}^{i_l^{\prime }})})^c$ and in particular, $%
Q_j^i\notin \cup _{l=1}^m\overline{W^s(Q_{j_l^{\prime }}^{i_l^{\prime }}).}$

\begin{proposition}
\label{propsab} As above, let $\cup _{l=1}^m\overline{W^u(Q_{j_l}^{i_l})}$
and $\cup _{l=1}^{m^{\prime }}\overline{W^s(Q_{j_l^{\prime }}^{i_l^{\prime
}})}$ be connected subsets, for integers $m,m^{\prime }\geq 1$ and sequences 
$1\leq i_l\leq s^{*}<i_l^{\prime }\leq n,$ $1\leq j_l,j_l^{\prime }\leq q.$

In case 
$$
({\rm Filled}(\cup _{l=1}^m\overline{W^u(Q_{j_l}^{i_l})}))^c\supset \cup
_{l=1}^{m^{\prime }}W^s(Q_{j_l^{\prime }}^{i_l^{\prime }})
$$
\begin{center}
and
\end{center}
$$
({\rm Filled}(\cup _{l=1}^{m^{\prime }}\overline{W^s(Q_{j_l^{\prime
}}^{i_l^{\prime }})}))^c\supset \cup _{l=1}^mW^u(Q_{j_l}^{i_l})), 
$$
which means that 
$$
{\rm Interior}\left( {\rm Filled}(\cup _{l=1}^{m^{\prime }}\overline{W^s(Q_{j_l^{\prime
}}^{i_l^{\prime }})})\right) \cap 
{\rm Interior}\left( {\rm Filled}(\cup _{l=1}^m\overline{W^u(Q_{j_l}^{i_l})})\right) =\emptyset,
$$
then the sum of topological indices of $f^q$ at all points in

${\rm Fix}(\widehat{f}^q(\bullet )-(0,p))\cap \left( {\rm Filled}(\cup _{l=1}^m 
\overline{W^u(Q_{j_l}^{i_l})})\cup {\rm Filled}(\cup _{l=1}^{m^{\prime }} 
\overline{W^s(Q_{j_l^{\prime }}^{i_l^{\prime }})})\right) $ is $2.$
\end{proposition}

{\it Proof. }

By contradiction, assume that the above sum of topological indices is less than $2.$
As Lemma \ref{gerale} implies that the sum of topological indices of $f^q$ at all points in
${\rm Fix}(\widehat{f}^q(\bullet )-(0,p))\cap \left( {\rm Filled}(\cup _{l=1}^m 
\overline{W^u(Q_{j_l}^{i_l})})\right)$ is equal to $1$ and 
the sum of topological indices of $f^q$ 
at all points in 
${\rm Fix}(\widehat{f}^q(\bullet )-(0,p))\cap \left( {\rm Filled}(\cup _{l=1}^{m^{\prime }} 
\overline{W^s(Q_{j_l^{\prime }}^{i_l^{\prime }})})\right)$ is also equal to $1,$ 
then there exists at least one point $w\in
{\rm Fix}(\widehat{f}^q(\bullet )-(0,p)),$ of topological index $1,$ belonging to both $\cup
_{l=1}^m\overline{W^u(Q_{j_l}^{i_l})}$ and $\cup _{l=1}^{m^{\prime }} 
\overline{W^s(Q_{j_l^{\prime }}^{i_l^{\prime }})}.$ As $f\in \chi ^r({\rm T^2%
}),$ $w$ must be an orientation reversing saddle, because it neither can be
a sink or a source. As $w\in \overline{W^u(Q_j^i)}\cap \overline{%
W^s(Q_{j^{\prime }}^{i^{\prime }})},$ for some $1\leq i\leq s^{*}<i^{\prime
}\leq n$ and $1\leq j,j^{\prime }\leq q,$ Theorem \ref{comeconovo} implies
that $W^u(Q_j^i)$ has a transversal intersection with $W^s(f^{j_1}(w))$ and $%
W^u(f^{j_2}(w))$ has a transversal intersection with $W^s(Q_{j^{\prime
}}^{i^{\prime }})$ for some $1\leq j_1,j_2\leq q.$ So, the $\lambda$-lemma
implies that $W^u(Q_j^i)$ accumulates on $W^u(Q_{j^{\prime \prime
}}^{i^{\prime }}),$ for some $1\leq j^{\prime \prime }\leq q,$ a
contradiction because $W^u(Q_j^i)$ is bounded when lifted to the plane, and $%
W^u(Q_{j^{\prime \prime }}^{i^{\prime }})$ is not. $\Box $

\vskip0.2truecm

Now, let us a consider a connected component of $\theta ,$ denoted $K,$
which can be decomposed as $K=\left( K_1^u\cup K_2^u\cup ...\cup
K_{n_u}^u\right) \cup \left( K_1^s\cup K_2^s\cup ...\cup K_{n_s}^s\right) ,$
where each $K_i^u$ is connected and given by the union of finitely many sets
of the form $\overline{W^u(Q_l^i)},$ $1\leq i\leq s^{*}$ and $1\leq l\leq q,$
and analogously, each $K_i^s$ is connected and given by the union of
finitely many sets of the form $\overline{W^s(Q_{l^{\prime }}^{i^{\prime }})}%
,$ $s^{*}<i^{\prime }\leq n$ and $1\leq l^{\prime }\leq q.$ The unstable 
 union is pairwise disjoint (similarly for the stable), that is,  
$K_i^u\cap K_j^u=\emptyset$ and $K_i^s\cap K_j^s=\emptyset,$ for 
$i\neq j.$ Clearly, either the
unstable or the stable part of $K$ could be empty.

Consider ${\rm Filled}(K_i^u)$ and ${\rm Filled}(K_j^u),$ for $i\neq j.$ The
possibilities are:

\begin{itemize}
\item  ${\rm Filled}(K_i^u)\subset {\rm Filled}(K_j^u);$

\item  ${\rm Filled}(K_i^u)\supset {\rm Filled}(K_j^u);$

\item  ${\rm Filled}(K_i^u)\cap {\rm Filled}(K_j^u)=\emptyset ;$
\end{itemize}

Analogously, for the stable components.

So, the set ${\rm Filled}(K_1^u)\cup ...\cup {\rm Filled}(K_{n_u}^u)$ might
not be a disjoint union anymore. But it becomes so, after omitting some
elements in the union, re-indexing the sets, and eventually decreasing the
value of $n_u$ (an analogous construction can be performed for the stable union).
On the other hand, any set of the form $W^s(Q_{l^{\prime
}}^{i^{\prime }}), $ contained in the stable union $K_1^s\cup K_2^s\cup
...\cup K_{n_s}^s$ is either contained in ${\rm Filled}(K_i^u)$ for some $%
1\leq i\leq n_u,$ or it belongs to $\left( {\rm Filled}(K_1^u)\cup ...\cup 
{\rm Filled}(K_{n_u}^u)\right) ^c,$ see Proposition \ref{propsab0}.

So we can remove from $K_1^s\cup K_2^s\cup ...\cup K_{n_s}^s$ all sets of
the form $\overline{W^s(Q_{l^{\prime }}^{i^{\prime }})}$ that are contained
in ${\rm Filled}(K_i^u),$ for some $1\leq i\leq n_u.$

Thus, there is a new stable union, still pairwise disjoint, denoted as $%
K_1^{*s}\cup K_2^{*s}\cup ...\cup K_{n_s^{\prime }}^{*s},$ such that each $%
K_i^{*s}$ is still connected and given by the union of finitely many sets of
the form $\overline{W^s(Q_{l^{\prime }}^{i^{\prime }})}$ $(s^{*}<i^{\prime
}\leq n$ and $1\leq l^{\prime }\leq q),$ satisfying the additional property
that each $\overline{W^s(Q_{l^{\prime }}^{i^{\prime }})}$ is contained in
the closure of $\left( {\rm Filled}(K_1^u)\cup ...\cup {\rm Filled}%
(K_{n_u}^u)\right) ^c. $ It is still possible that ${\rm Filled}%
(K_j^{*s})\supset {\rm Filled}(K_i^u)$ for indices $1\leq j\leq n_s^{\prime
} $ and $1\leq i\leq n_u.$ In this case, we just omit $K_i^u$ from the
unstable union. After all these modifications, we end up with the following
set: 
\begin{equation}
\label{Kchapeu}\widehat{K}= [{\rm Filled}(K_1^u)\cup ...\cup {\rm Filled}%
(K_{n_u}^u)] \cup [{\rm Filled}(K_1^{*s})\cup ...\cup {\rm Filled}%
(K_{n_s^{\prime }}^{*s})] , 
\end{equation}
satisfying the conditions below:

\begin{enumerate}
\item  $\widehat{K}$ is closed, connected and contains $K;$

\item  its unstable components are still pairwise disjoint, as are the
stable;

\item  ${\rm interior}({\rm Filled}(K_i^u))$ does not intersect ${\rm %
interior}({\rm Filled}(K_j^{*s}));$

\item $\widehat{K}^c$ has a fully essential component;
\end{enumerate}

Any connected component of $\theta $ (see expression (\ref{defteta})), is
either contained in $\widehat{K}$ or is disjoint from it, because $\partial 
\widehat{K}\subset K\subset \theta $.

Now let $M$ be a connected component of $\theta $ which is disjoint from $
\widehat{K}.$ We want to show that the set $\widehat{M},$ constructed from $%
M $ exactly in the same way as $\widehat{K}$ was constructed from $K,$
satisfies:\ either $\widehat{M}\cap \widehat{K}=\emptyset $ or $\widehat{M}%
\supset \widehat{K.}$ But this is easy: if $\widehat{M}\cap \widehat{K}\neq
\emptyset ,$ as we assumed that $M\cap \widehat{K}=\emptyset ,$ from the
fact that $\partial \widehat{M}\subset M,$ we get that $\partial \widehat{M}$
avoids $\widehat{K}.$ So, as $\widehat{K}$ is connected, $\widehat{M}\supset 
\widehat{K}.$

Therefore, we can find finitely many disjoint connected closed sets of the
form $\widehat{K}$ as above (see expression (\ref{Kchapeu})), such that $%
\theta $ is contained in their union and 
\begin{equation}
\label{expfina}\sum_{z\in {\rm Fix}(\widehat{f}^q(\bullet )-(0,p))\cap \widehat{K}%
}{\rm ind}(f^q\mid _z)=n_u+n_s^{\prime }. 
\end{equation}
As each $\widehat{K}$ contains some connected component of $\theta ,$ $%
n_u+n_s^{\prime }\geq 1.$

In order to show that expression (\ref{expfina}) holds, note that Lemma \ref
{gerale} implies that 
$$
\sum_{z\in {\rm Fix}(\widehat{f}^q(\bullet )-(0,p))\cap 
{\rm Filled}(K_i^u)}{\rm ind}(f^q\mid _z)=1 
$$
$$
\text{ and}
$$ 
$$
\sum_{z\in {\rm Fix}(\widehat{f}
^q(\bullet )-(0,p))\cap {\rm Filled}(K_i^{*s})}{\rm ind}(f^q\mid _z)=1.
$$

If the
sum of topological indices in expression (\ref{expfina}) were strictly less than $%
n_u+n_s^{\prime },$ as unstable unions are pairwise disjoint (as are the
stable ones), a contradiction with Proposition \ref{propsab} would have been
found.

So, from the fact that $\theta $ contains all $q$-periodic saddles of
vertical rotation number $p/q$ and topological index $-1,$ this contradicts expression
(\ref{lefschetz}) and proves the theorem. $\Box $

\vskip0.2truecm

\subsection{Proof of Theorem 4 and Corollary 2}

From the introduction: Let $k\neq 0$ be an integer and $\widehat{f}\in O{\rm _{t,k}^r}({\rm %
T^1\times I\negthinspace  R})\cap \chi ^r({\rm T^2})$ or 
$\widehat{f}\in O{\rm _{t,k,Leb}^r}({\rm T^1\times I%
\negthinspace  R}) \cap \chi _{{\rm Leb}}^r({\rm T^2}),$ 
the open and dense sets from Theorems \ref{conjrot} and 
\ref{conjrotarea} intersected with the generic sets from Theorem \ref{comeconovo}. 
So $\rho _V(\widehat{f})=[r/s,p/q]$ for
rational numbers $r/s\leq $ $p/q$ and $\rho _V$ is locally constant in a
neighborhood of $\widehat{f}.$ Assume that $r/s<p/q.$ Then, $
\widehat{f}^q(\bullet )-(0,p)$ has a free homotopically non-trivial simple
closed curve $\gamma _{p/q}\subset {\rm T^1\times I\negthinspace  R}$, such
that: 
$$
\widehat{f}^q(\gamma _{p/q})-(0,p)\subset \gamma _{p/q}^{-} 
$$

Something which implies the existence of an attractor-repeller pair for $
\widehat{f}^q(\bullet )-(0,p).$ The attractor $A_{p/q}$ is contained in $%
\gamma _{p/q}^{-}$ and the repeller $R_{p/q}$ is contained in $\gamma
_{p/q}^{+}.$ The statements proved here are:

\begin{enumerate}
\item  Under the previous hypotheses, there exists a hyperbolic periodic
saddle $z_{p/q}\in {\rm T^2}$ of vertical rotation number $p/q$ (and
period $nq$ for some $n\geq 1$) such that if $\widehat{z}_{p/q}\in {\rm %
T^1\times I\negthinspace  R}$ is any lift of $z_{p/q}$ to the annulus, then $%
W^u(\widehat{z}_{p/q})$ is bounded from above as a subset of the annulus and
unbounded from below, $W^s(\widehat{z}_{p/q})$ is unbounded from above and
bounded from below. Moreover, $W^u(\widehat{z}_{p/q})$ has a transversal
intersection with $W^s(\widehat{z}_{p/q}-(0,1))=W^s(\widehat{z}%
_{p/q})-(0,1), $ and if $\widehat{z}_{p/q}$ and $\widehat{z}_{p/q}-(0,1)$
are both above $\gamma _{p/q},$ and $\widehat{z}_{p/q}-(0,2)$ is below, then 
$$
A_{p/q}\subset \overline{W^u(\widehat{z}_{p/q})}\cup (\overline{W^u(\widehat{%
z}_{p/q})}^{b.above}), 
$$
where the last set is the union of all (open) connected components of $( 
\overline{W^u(\widehat{z}_{p/q})})^c$ which are bounded from above.

Moreover, $A_{p/q}\supseteq \left[ \overline{W^u(\widehat{z}_{p/q})}\cup ( 
\overline{W^u(\widehat{z}_{p/q})}^{b.above})\right] -(0,2).$

Similarly, if $\widehat{z}_{p/q}$ and $\widehat{z}_{p/q}+(0,1)$ are both
below $\gamma _{p/q},$ and $\widehat{z}_{p/q}+(0,2)$ is above, then 
$$
R_{p/q}\subset \overline{W^s(\widehat{z}_{p/q})}\cup (\overline{W^s(\widehat{%
z}_{p/q})}^{b.below}), 
$$
where the last set is the union of all connected components of $(\overline{%
W^s(\widehat{z}_{p/q})})^c$ which are bounded from below.

And analogously, $R_{p/q}\supseteq \left[ \overline{W^s(\widehat{z}_{p/q})}%
\cup (\overline{W^s(\widehat{z}_{p/q})}^{b.below})\right] +(0,2).$

\item  If $f$ is transitive, then the following improvement holds:

a) assuming that $\widehat{z}_{p/q}$ and $\widehat{z}_{p/q}-(0,1)$ are both
above $\gamma _{p/q},$ and $\widehat{z}_{p/q}-(0,2)$ is below, then $A_{p/q}$
is contained in $\overline{W^u(\widehat{z}_{p/q})}$ and contains $\overline{%
W^u(\widehat{z}_{p/q})}-(0,2);$

b) similarly, if $\widehat{z}_{p/q}$ and $\widehat{z}_{p/q}+(0,1)$ are both
below $\gamma _{p/q},$ and $\widehat{z}_{p/q}+(0,2)$ is above, then $R_{p/q}$
is contained in $\overline{W^s(\widehat{z}_{p/q})}$ and contains $\overline{%
W^s(\widehat{z}_{p/q})}+(0,2);$

c) Moreover, both $\overline{W^u(\widehat{z}_{p/q})}$ and $\overline{W^s(
\widehat{z}_{p/q})}$ have no interior points and their complements are
connected. In the torus,  
$\overline{W^u(z_{p/q})}=\overline{W^s(z_{p/q})}={\rm T^2}.$ 
\end{enumerate}

{\it Proof of 1).}

As $\rho _V(\widehat{f})=[r/s,p/q]$ is locally constant in a $C^r$%
-neighborhood of $\widehat{f}\in {\rm Diff_{t,k}^r(T^1\times I\negthinspace  %
R),}$ Theorem \ref{rational1} implies that $f$ has $nq$-periodic points of
vertical rotation number $p/q$ for some integer $n\geq 1.$ So, Theorem \ref
{gerall} implies the existence of a hyperbolic $f$-periodic saddle point $%
z_{p/q}$ of vertical rotation number $p/q$ (and period $nq),$ with a
partial mesh. Thus, $\widehat{f}^q(\bullet )-(0,p)$ has a $n$-periodic
hyperbolic saddle point $\widehat{z}_{p/q}\in {\rm T^1\times I\negthinspace  %
R,}$ for which, as $\rho _V^{+}(\widehat{f})=p/q,$ Corollary 1 of \cite
{braulio} implies that $W^u(\widehat{z}_{p/q})$ is bounded from above as a
subset of the annulus and $W^s(\widehat{z}_{p/q})$ is bounded from below.
Clearly all integer vertical translates of $\widehat{z}_{p/q}$ are $n$%
-periodic for $\widehat{f}^q(\bullet )-(0,p)$ and satisfy the above
properties.

\begin{description}
\item[Fact 1.]  {\it The manifolds  
 $W^u(\widehat{z}_{p/q})$ and $W^s(\widehat{z}_{p/q})-(0,1)$} intersect transversely.
\end{description}

{\it Proof.}

In order to see this, note that as $z_{p/q}\in {\rm T^2}$ has a partial
mesh, and $\rho _V^{+}(\widehat{f})=p/q,$ one of the following possibilities
holds (maybe both):

\begin{enumerate}
\item  $W^u(\widehat{z}_{p/q})\cup W^s(\widehat{z}_{p/q})$ contains a
homotopically non-trivial closed curve $\kappa $ in the annulus;

\item  $W^u(\widehat{z}_{p/q})$ intersects $W^s(\widehat{z}_{p/q})-(0,l)$
for integers $l>0;$
\end{enumerate}

These are the only possibilities, because if $W^u(\widehat{z}_{p/q})$
intersected $W^s(\widehat{z}_{p/q})+(0,i)$ for some integer $i>0,$ then $%
\rho _V^{+}(\widehat{f})$ would be larger than $p/q,$ because such an
intersection would imply the existence of a rotational horseshoe for $
\widehat{f}^q(\bullet )-(0,p)$ and this horseshoe would contain periodic
points with positive vertical rotation number with respect to the map $
\widehat{f}^q(\bullet )-(0,p)$, a contradiction.

In the first case above, if $\kappa $ and $\kappa -(0,1)$ intersect, as $f$
is Kupka-Smale, then $\kappa \cap $$W^u(\widehat{z}_{p/q})$ must have a
transversal intersection with $(\kappa -(0,1))\cap (W^s(\widehat{z}%
_{p/q})-(0,1))$ and we are done. If they do not intersect, as $\rho _V^{-}( 
\widehat{f})=r/s<\rho _V^{+}(\widehat{f})=p/q,$ then 
$$
\left( \widehat{f}^{nq}(\bullet )-(0,np)\right) ^m(\kappa )\cap (\kappa
-(0,1))\neq \emptyset 
$$
for some integer $m>0,$ and we get the same conclusion as before.

So let us assume by contradiction, that $W^u(\widehat{z}_{p/q})$ intersects $%
W^s(\widehat{z}_{p/q})-(0,l_0)$ for some integer $l_0>1$ and $W^u(\widehat{z}%
_{p/q})$ does not intersect $W^s(\widehat{z}_{p/q})-(0,i)$ for $%
i=1,2,...,l_0-1.$

First recall that from the choice of $f,$ it is homotopic to 
$$
(x,y)\mapsto \left( 
\begin{array}{cc}
1 & k \\ 
0 & 1 
\end{array}
\right) \left( 
\begin{array}{c}
x \\ 
y 
\end{array}
\right) (\rm{mod\ 1})^2. 
$$

From our assumption, if $\widetilde{z}_{p/q}$ is a lift of $z_{p/q}$ to the
plane, then there exists an integer $a$ and a simple arc $\widetilde{\zeta }%
, $ starting at $\widetilde{z}_{p/q}$ and ending at $\widetilde{z}%
_{p/q}-(a,l_0)$ made of two connected pieces: the first one, $\widetilde{%
\zeta }_1,$ is contained in $W^u(\widetilde{z}_{p/q}),$ starts at $
\widetilde{z}_{p/q}$ and ends at a point $\widetilde{w}\in W^u(\widetilde{z}%
_{p/q})\cap (W^s(\widetilde{z}_{p/q})-(a,l_0)),$ and the second one, $
\widetilde{\zeta }_2,$ is contained in $W^s(\widetilde{z}_{p/q})-(a,l_0),$
starts at $\widetilde{w}$ and ends at $\widetilde{z}_{p/q}-(a,l_0).$ Now, fixed 
some $\widetilde{f}:{\rm I\negthinspace R^2}\rightarrow {\rm I\negthinspace R^2,}$
 lift of $\widehat{f}$ to the plane,
consider the image of $\widetilde{\zeta }$ under $\widetilde{f}^{nq}(\bullet
)-(s^{\prime },np),$ where the integer $s^{\prime }$ was chosen so that $
\widetilde{z}_{p/q}$ is fixed: $\widetilde{f}^{nq}(\widetilde{\zeta }%
_1)-(s^{\prime },np)$ contains $\widetilde{\zeta }_1$ and $\widetilde{f}%
^{nq}(\widetilde{\zeta }_2)-(s^{\prime },np)$ is contained in $\widetilde{%
\zeta }_2-(nqkl_0,0).$ So, $(\widetilde{f}^{nq}(\widetilde{\zeta }%
_1)-(s^{\prime },np))\cup \widetilde{\zeta }_2\cup (\widetilde{\zeta }%
_2-(nqkl_0,0))$ contains an arc connecting $\widetilde{z}_{p/q}-(a,l_0)$ to $
\widetilde{z}_{p/q}-(a+nqkl_0,l_0).$ Therefore, $W^u(\widehat{z}_{p/q})\cup
(W^s(\widehat{z}_{p/q})-(0,l_0))$ contains a homotopically non-trivial
simple closed curve $\mu $ (every homotopically non-trivial closed curve in
the annulus contains a homotopically non-trivial SIMPLE closed curve). Let
us consider the point $\widehat{z}_{p/q}-(0,1).$ As there are no saddle
connections, it does not belong to $\mu .$ If it is below $\mu ,$ then as
its stable manifold is unbounded from above, it must intersect $\mu ,$ and
so $W^u(\widehat{z}_{p/q})$ intersects $W^s(\widehat{z}_{p/q})-(0,1),$ a
contradiction. And if $\widehat{z}_{p/q}-(0,1)$ is above $\mu ,$ then as its
unstable manifold is unbounded from below, it must intersect $\mu ,$ and so $%
W^u(\widehat{z}_{p/q})-(0,1)$ intersects $W^s(\widehat{z}_{p/q})-(0,l_0),$
which implies that $W^u(\widehat{z}_{p/q})$ intersects $W^s(\widehat{z}%
_{p/q})-(0,l_0-1),$ a contradiction with the choice of $l_0.$ $\Box $

\vskip0.2truecm

As $\widehat{f}^{nq}(\bullet )-(0,np)$ fixes all integer vertical translates
of $\widehat{z}_{p/q},$ and maps $\gamma _{p/q}$ into $\gamma _{p/q}^{-},$ the
lower connected component of $\gamma _{p/q}^c,$ we get that for all integers 
$i,$ $\widehat{z}_{p/q}+(0,i)$ does not belong to $\gamma _{p/q}.$ Suppose
for some integer $i,$ $\widehat{z}_{p/q}+(0,i)$ is below $\gamma _{p/q}.$
Then $W^u(\widehat{z}_{p/q})+(0,i)$ is also contained in $\gamma _{p/q}^{-}.$
Otherwise, it would intersect $\gamma _{p/q}$ in a point $\widehat{z}%
^{\prime }.$ The negative orbit of $\widehat{z}^{\prime }$ under $\widehat{f}%
^{nq}(\bullet )-(0,np)$ converges to $\widehat{z}_{p/q}+(0,i),$ and is
always above $\gamma _{p/q},$ so $\widehat{z}_{p/q}+(0,i)$ belongs to $%
\gamma _{p/q}^{+},$ a contradiction with its choice. Similarly, if for some
integer $i,$ $\widehat{z}_{p/q}+(0,i)$ is above $\gamma _{p/q},$ then $W^s( 
\widehat{z}_{p/q})+(0,i)$ is also contained in $\gamma _{p/q}^{+}.$

The following two items describe how the integer vertical translates of $
\widehat{z}_{p/q}$ behave when compared to $\gamma _{p/q}.$

\begin{itemize}
\item  if for some integer $i,$ $\widehat{z}_{p/q}+(0,i)$ is above $\gamma
_{p/q},$ then for all integers $j>i,$ $\widehat{z}_{p/q}+(0,j)$ is also
above $\gamma _{p/q}.$
\end{itemize}

To see this, recall that for all integers $l,$ $W^u(\widehat{z}%
_{p/q})+(0,l)$ intersects $W^s(\widehat{z}_{p/q})+(0,l-1)$ transversely. So
the $\lambda $-lemma implies that $W^u(\widehat{z}_{p/q})+(0,l)$ intersects $%
W^s(\widehat{z}_{p/q})+(0,l^{\prime })$ for all integers $l^{\prime }<l.$
Thus, $W^u(\widehat{z}_{p/q}+(0,l))=W^u(\widehat{z}_{p/q})+(0,l)$
accumulates on $W^u(\widehat{z}_{p/q})+(0,l^{\prime })$ for all $l^{\prime
}<l.$ Therefore, if $\widehat{z}_{p/q}+(0,i)$ is above $\gamma _{p/q}$ and $%
j>i,$ $W^u(\widehat{z}_{p/q})+(0,j)$ accumulates on $W^u(\widehat{z}%
_{p/q})+(0,i),$ something that implies that $W^u(\widehat{z}_{p/q})+(0,j)$
intersects $\gamma _{p/q}^{+}.$ So $\widehat{z}_{p/q}+(0,j)$ is not below $%
\gamma _{p/q},$ and we are done.

\begin{itemize}
\item  similarly, if for some integer $i,$ $\widehat{z}_{p/q}+(0,i)$ is
below $\gamma _{p/q},$ then for all integers $j<i,$ $\widehat{z}_{p/q}+(0,j)$
is also below $\gamma _{p/q}.$
\end{itemize}

The proof is as above.

Summarizing, when considering integer vertical translates of $\widehat{z}%
_{p/q},$ there exists an integer $\overline{i}$ such that for all $i\leq 
\overline{i},$ $\widehat{z}_{p/q}+(0,i)$ is below $\gamma _{p/q}$ and for
all $i>\overline{i},$ $\widehat{z}_{p/q}+(0,i)$ is above $\gamma _{p/q}.$

By re-indexing integer vertical translates of $\widehat{z}_{p/q},$ we can
assume that $\widehat{z}_{p/q}$ and $\widehat{z}_{p/q}-(0,1)$ are both above 
$\gamma _{p/q}$ and $\widehat{z}_{p/q}-(0,2)$ is below. From the proof of
Fact 1 stated above, $W^u(\widehat{z}_{p/q})\cup (W^s(\widehat{z}_{p/q})-(0,1))$
contains a homotopically non-trivial simple closed curve $\mu .$

For some large integer $m^{*}>0,$ $\left( \widehat{f}^{nq}(\bullet
)-(0,np)\right) ^{-m^{*}}(\mu )$ is above $\gamma _{p/q}.$

To see this, recall that $\mu $ is contained in $W^u(\widehat{z}%
_{p/q})\cup (W^s(\widehat{z}_{p/q})-(0,1)).$ So, if $m>0$ is large enough, 
$$
\left( \widehat{f}^{nq}(\bullet )-(0,np)\right) ^{-m}(\mu )\cap W^u(\widehat{%
z}_{p/q}) 
$$
is contained in a small neighborhood of $\widehat{z}_{p/q},$ and thus, it is
contained in $\gamma _{p/q}^+.$ As the whole $W^s(\widehat{z}_{p/q})-(0,1)$
is above $\gamma _{p/q},$ the result follows.

But now, if we denote $\left( \widehat{f}^{nq}(\bullet )-(0,np)\right)
^{-m^{*}}(\mu )$ as $\mu _{*},$ then $\mu _{*}$ is a homotopically
non-trivial simple closed curve such that $\mu _{*}^{-}\supset \gamma
_{p/q}. $

So, for all integers $j\geq 0,$ 
\begin{equation}
\label{fimitem1} 
\begin{array}{c}
A_{p/q}:=\cap _{i=0}^{+\infty }\left( 
\widehat{f}^q(\bullet )-(0,p)\right) ^i(\gamma _{p/q}^{-})\subset \left( 
\widehat{f}^{nq}(\bullet )-(0,np)\right) ^j(\gamma _{p/q}^{-})\subset \\ 
\subset \left( \widehat{f}^{nq}(\bullet )-(0,np)\right) ^j(\mu _{*}^{-}). 
\end{array}
\end{equation}
Moreover, as $\widehat{z}_{p/q}-(0,2)$ belongs to $\gamma _{p/q}^{-},$ $W^u( 
\widehat{z}_{p/q})-(0,2)$ is also contained in $\gamma _{p/q}^{-}.$ As $W^u( 
\widehat{z}_{p/q})-(0,2)$ is $\left( \widehat{f}^{nq}(\bullet
)-(0,np)\right) $-invariant, its closure is contained in the closed set $%
A_{p/q}.$ And clearly, from the definition of $A_{p/q},$ as it contains $
\overline{W^u(\widehat{z}_{p/q})}-(0,2),$ then $A_{p/q}$ also contains $( 
\overline{W^u(\widehat{z}_{p/q})}^{b.above})-(0,2),$ where $\overline{W^u( 
\widehat{z}_{p/q})}^{b.above}$ was previously defined.

In order to finish the proof of 1, we are left to show that expression (\ref
{fimitem1}) implies that $\overline{W^u(\widehat{z}_{p/q})}\cup (\overline{%
W^u(\widehat{z}_{p/q})}^{b.above})$ contains $A_{p/q}.$

It clearly implies that 
$$
\cap _{i=0}^{+\infty }\left( \widehat{f}^{nq}(\bullet )-(0,np)\right) ^i(%
{\rm closure}(\mu _{*}^{-}))\supset A_{p/q}, 
$$
so let us show the following:

\begin{description}
\item[Fact 2.] {\it The set $\theta ^{*}:=\cap _{i=0}^{+\infty }\left( \widehat{f%
}^{nq}(\bullet )-(0,np)\right) ^i({\rm closure}(\mu _{*}^{-}))$ is contained
in $\overline{W^u(\widehat{z}_{p/q})}$ $\cup (\overline{W^u(\widehat{z}%
_{p/q})}^{b.above}).$}
\end{description}

{\it Proof.}

First, recall that $\mu _{*}\subset W^u(\widehat{z}_{p/q})\cup (W^s( 
\widehat{z}_{p/q})-(0,1)),$ and pick some $\widehat{w}\in \theta ^{*}.$
Suppose it does not belong to $\overline{W^u(\widehat{z}_{p/q})}.$ Then it
belongs to some connected component of $(\overline{W^u(\widehat{z}_{p/q})}%
)^c.$ If $\widehat{w}\notin \overline{W^u(\widehat{z}_{p/q})}^{b.above},$
then $\widehat{w}$ belongs to a connected component of $(\overline{W^u( 
\widehat{z}_{p/q})})^c$ which is not bounded from above. As $\overline{W^u( 
\widehat{z}_{p/q})}$ is itself bounded from above, there is only one
connected component of $(\overline{W^u(\widehat{z}_{p/q})})^c$ which is
unbounded from above, denoted $Unb.$ Let $a$ be a real number such that $%
{\rm T^1}\times \{a\}\subset Unb$ and ${\rm T^1}\times \{a\}$ avoids $\left( \widehat{f}%
^{nq}(\bullet )-(0,np)\right) ^j(\mu _{*})$ for all integers $j\geq 0.$ Let $%
\psi \subset Unb$ be a simple arc connecting $\widehat{w}$ to some point
above ${\rm T^1}\times \{a\}.$

As $\widehat{z}_{p/q}-(0,1)\in \overline{W^u(\widehat{z}_{p/q})},$ if $j>0$
is large enough, then 
\begin{equation}
\label{consertouhora}\left[ \left( \widehat{f}^{nq}(\bullet )-(0,np)\right)
^j(\mu _{*})\cap \left( W^s(\widehat{z}_{p/q})-(0,1)\right) \right] \cap
\psi =\emptyset . 
\end{equation}
And as $\widehat{w}\in \theta ^{*}$ and the other endpoint of $\psi $ is
above ${\rm T^1}\times \{a\},$ we get that $\psi $ intersects $\left( \widehat{f}%
^{nq}(\bullet )-(0,np)\right) ^j(\mu _{*})$ for all $j\geq 0.$ So, if $j>0$
is large enough, expression (\ref{consertouhora}) implies that $\left( 
\widehat{f}^{nq}(\bullet )-(0,np)\right) ^j(\mu _{*})\cap \psi $ belongs to $%
W^u(\widehat{z}_{p/q}),$ a contradiction with the choice of $\psi$ that
finishes the proof. 
$\Box $

\vskip0.2truecm

Item 1) is thus proved, in case of the attractor $A_{p/q}.$ The proof for
the repeller $R_{p/q}$ is analogous because $R_{p/q}$ is an attractor for $%
\left( \widehat{f}^q(\bullet )-(0,p)\right) ^{-1}.$ $\Box $

\vskip0.2truecm

{\it Proof of 2).}

The first thing to prove is, in case $f:{\rm T^2\rightarrow T^2}$ is
transitive, that both $W^u(z_{p/q})$ and $W^s(z_{p/q})$ are dense in ${\rm T^2.}$
This is already known when $f$ preserves area, see for instance \cite
{malha}. Let us prove it in general. Suppose $f\in {\rm Diff_{t,k}^r(T^2)%
\cap \chi ^r(T^2)}$ and, for instance, $\overline{W^u(z_{p/q})}$ is not the
whole torus. In this case, as $z_{p/q}$ has a partial mesh, any connected
component of its complement is a topological open disk. Let us consider such
a disk $D,$ which can be either periodic, or wandering. If $f$ is
transitive, then $D$ cannot be wandering. So, it is $n$-periodic, for some
integer $n>0.$ Moreover, $D$ cannot have homotopically bounded diameter
(that is, any connected components of $p^{-1}(D)$ has the same bounded
diameter), because $\overline{D}\cup f(\overline{D})\cup ...\cup f^{n-1}(
\overline{D})={\rm T^2,}$ and if $D$ has homotopically bounded diameter,
then all points in ${\rm T^2}$ would have the same vertical rotation
number (because points in the boundary of a homotopically bounded periodic
disk have the same vertical rotation number as points in the disk), a
contradiction with the assumption that $\rho _V(\widehat{f})$ is a
non-degenerate interval. And if $D$ is homotopically unbounded, Proposition
24 of \cite{euandres} shows that there exists a homotopically bounded open
disk $D^{*}\subset D,$ such that $f^n(D^{*})\subset D^{*},$ where $n$ is the
period of $D.$ And this again contradicts the transitivity of $f$ and the
assumption that $\rho _V(\widehat{f})$ is a non-degenerate interval.

Let $\widehat{z}_{p/q}\in {\rm T^1\times I\negthinspace  R}$ be any point in
the fiber of $z_{p/q}.$ If by contradiction, $(\overline{W^u(\widehat{z}%
_{p/q})})^c$ is not connected, then it has a connected component $B^{+}$
which is bounded from above. So the whole orbit of $B^{+}$ under $\widehat{f}%
^{nq}(\bullet )-(0,np)$ is bounded from above, because the boundary of any
iterate of $B^{+}$ is contained in $\overline{W^u(\widehat{z}_{p/q})},$
which is itself bounded from above. As $W^u(z_{p/q})$ is dense in ${\rm T^2}$
and $W^u(\widehat{z}_{p/q})$ is bounded from above and it has a transverse
intersection with $W^s(\widehat{z}_{p/q})-(0,1),$ we get that for any
sufficiently large integer $l>0,$ $W^u(\widehat{z}_{p/q})+(0,l)$ intersects $%
B^{+}.$ So, for all sufficiently large $i>0,$ $\left( \widehat{f}%
^{nq}(\bullet )-(0,np)\right) ^{-i}(B^{+})$ gets close to $\widehat{z}%
_{p/q}+(0,l)$ by less then one. As $l>0$ is arbitrarily large, this is a
contradiction. So $(\overline{W^u(\widehat{z}_{p/q})})^c$ is connected, and
an easy modification of this argument shows that $\overline{W^u(\widehat{z}%
_{p/q})}$ has no interior. Analogous results hold for the stable manifold.
This proves sub-items a), b) and c). $\Box $

\vskip0.2truecm

\vskip0.2truecm

\textit{Acknowledgments:} I would like to warmly thank the referee for 
his very careful reading of the paper. Particularly, for the suggestions 
concerning the proofs of Theorems 2 and 3. The clarity of the argument and
the readability were greatly improved. I am also very thankful for all his 
other comments and remarks. It is not usual to see such a careful work.

\end{document}